\input amstex
\input amsppt.sty

\input epsf
\epsfverbosetrue  \magnification
980\vsize=21 true cm \hsize=16.5 true cm \voffset=1.1 true cm
\pageno=1 \NoRunningHeads \TagsOnRight

\def\p{\partial}
\def\ve{\varepsilon}
\def\f{\frac}
\def\na{\nabla}
\def\la{\lambda}
\def\al{\alpha}
\def\t{\tilde}
\def\vp{\varphi}
\def\O{\Omega}
\def\th{\theta}
\def\g{\gamma}

\def\a{(b_0 q_0)^{\f{\g-3}{\g-1}}}
\def\b{(b_0 q_0)^{-\f{2}{\g-1}}}
\def\c{(b_0 q_0)^{-2}}

\def\ds{\displaystyle}

\topmatter
\topmatter \vskip 0.3 true cm \title{\bf On the global existence
 and stability of a three-dimensional supersonic conic shock wave}
\endtitle
\endtopmatter
\document

\vskip 0.2 true cm \footnote""{* Li Jun and Yin Huicheng were
supported by the NSFC (No.~10931007, No.~11025105, No.~11001122), by
the Doctoral Program Foundation of the Ministry of Education of China
(No.~20090091110005), and by the DFG via the Sino-German project
``Analysis of PDEs and application.'' This research was carried out
while Li Jun and Yin Huicheng were visiting the Mathematical Institute
of the University of G\"{o}ttingen.} \footnote""{** Ingo Witt was
partly supported by the DFG via the Sino-German project ``Analysis of
PDEs and application.'' } \vskip 0.3 true cm \centerline{Li,
Jun$^{1,*}$; \qquad Witt, Ingo$^{2,**}$;\qquad Yin,
Huicheng$^{1,*}$} \vskip 0.5 true cm {1. Department of Mathematics and
IMS, Nanjing University, Nanjing 210093, P.R.~China.}\vskip 0.2 true
cm {2. Mathematical Institute, University of G\"{o}ttingen,
Bunsenstr.~3-5, D-37073 G\"{o}ttingen, Germany.} \vskip 0.4 true cm

\centerline {\bf Abstract} We establish the global existence and
stability of a three-dimensional supersonic conic shock wave for a
perturbed steady supersonic flow past an infinitely long circular cone
with a sharp angle. The flow is described by a 3-D steady potential
equation, which is multi-dimensional, quasilinear, and hyperbolic with
respect to the supersonic direction. Making use of the geometric
properties of the pointed shock surface together with the
Rankine-Hugoniot conditions on the conic shock surface and the
boundary condition on the surface of the cone, we obtain a global
uniform weighted energy estimate for the nonlinear problem by finding
an appropriate multiplier and establishing a new Hardy-type inequality
on the shock surface. Based on this, we prove that a multi-dimensional
conic shock attached to the vertex of the cone exists globally when
the Mach number of the incoming supersonic flow is sufficiently
large. Moreover, the asymptotic behavior of the 3-D supersonic conic
shock solution, that is shown to approach the corresponding background
shock solution in the downstream domain for the uniform supersonic
constant flow past the sharp cone, is also explicitly given.

\vskip 0.2 true cm

\vskip 0.2 true cm

\vskip 0.2 true cm

\vskip 0.4 true cm \centerline{\bf \S1. Introduction} \vskip
0.2 true cm

In this paper, we are concerned with the multi-dimensional steady and
supersonic conic shock wave problem for a perturbed incoming
supersonic flow past an infinitely long circular cone. This problem is
fundamental in gas dynamics, for instance, for the supersonic flight
of projectiles and rockets. It is also one of the basic models for the
discussion of the theory of weak solutions to quasilinear hyperbolic
equations and systems in multi-dimensions (see [3], [22]--[23], [30],
[37]). Under suitable assumptions on the incoming supersonic flow with
a small spherically symmetric perturbation and a spherically symmetric
pointed body or artificial boundary conditions on the conic surface,
there is an extensive literature studying supersonic flow past a
pointed body (see [5], [7]--[10], [18], [21], [31], [35]--[36], and
the references therein). The first rigorous mathematical analysis was
given in [7] by Courant and Friedrichs, who proved that, for a uniform
supersonic flow $(0,0,q_0)$ with constant density $\rho_0>0$ which
approaches from minus infinity, when the flow hits the sharp circular
cone $\sqrt{x_1^2+x_2^2}=b_0x_3$, $b_0>0$, in direction of the
$x_3$-axis (see Figure~1 below), then there appears a supersonic conic
shock $\sqrt{x_1^2+x_2^2}=s_0x_3$, $s_0>b_0$, attached to the tip of
the cone provided that $b_0$ is less than some critical value $b^*>0$,
which is determined by the parameters of the incoming flow. When the
incoming supersonic flow is multi-dimensionally perturbed, the basic
problem of both mathematical and physical relevance that naturally
arises is whether such a conic shock is globally stable. Or else, do
there appear new shocks or other complicated singularities in the
downstream domain? Here, we will focus on this problem when the Mach
number of the incoming supersonic flow is appropriately large. It will
be shown that a global supersonic conic shock exists uniquely in the
whole space and that there are no other singularities between the
conic shock and the conic surface for a multi-dimensionally perturbed
supersonic polytropic gas past the sharp cone
$\sqrt{x_1^2+x_2^2}=b_0x_3$ (see Figure~2 below). This result agrees
with physical experiment and numerical simulations.
$$
\epsfysize=75mm \epsfbox{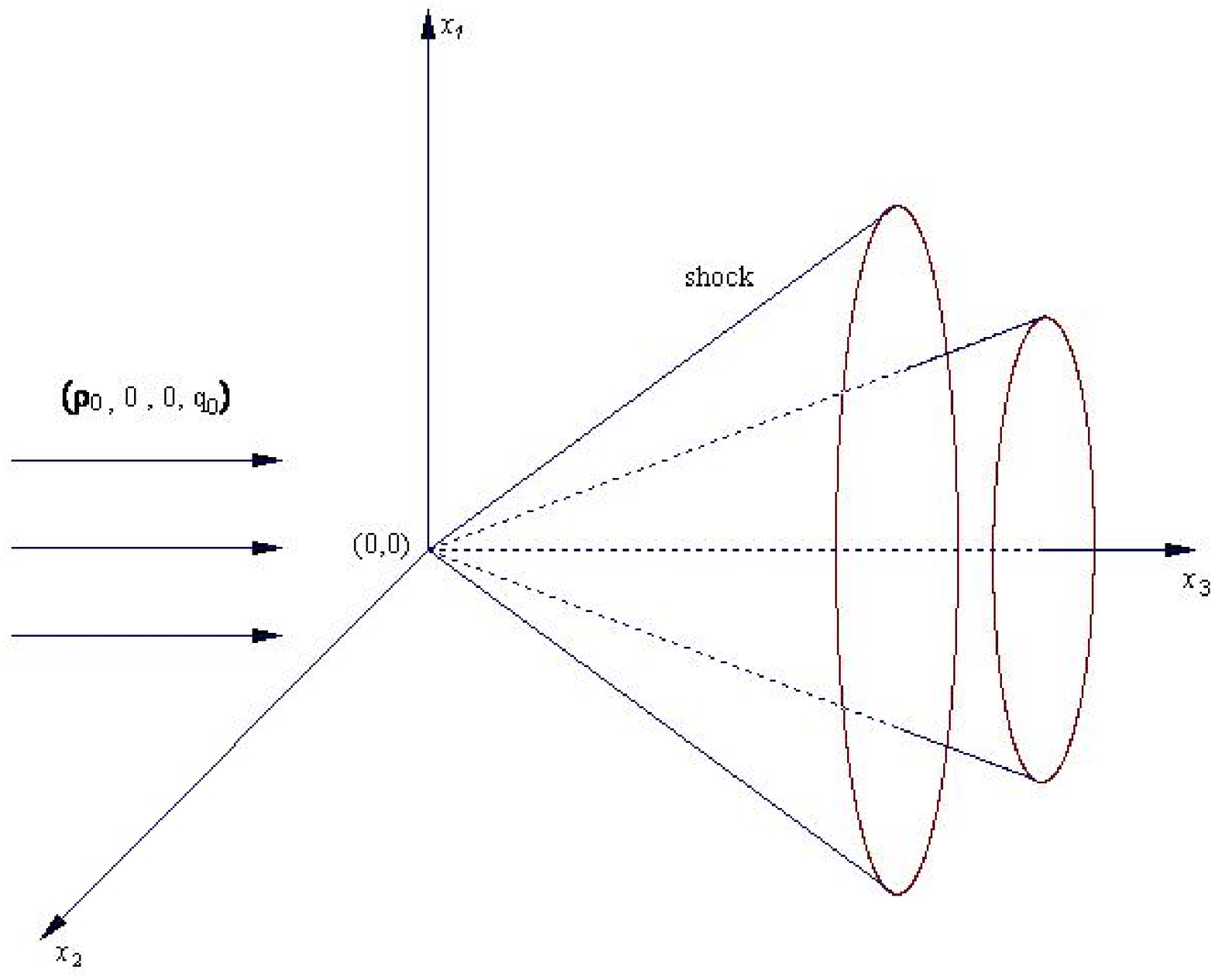}
$$
\centerline{\bf Figure 1. A uniform supersonic flow past a sharp
circular cone}
$$
\epsfysize=75mm \epsfbox{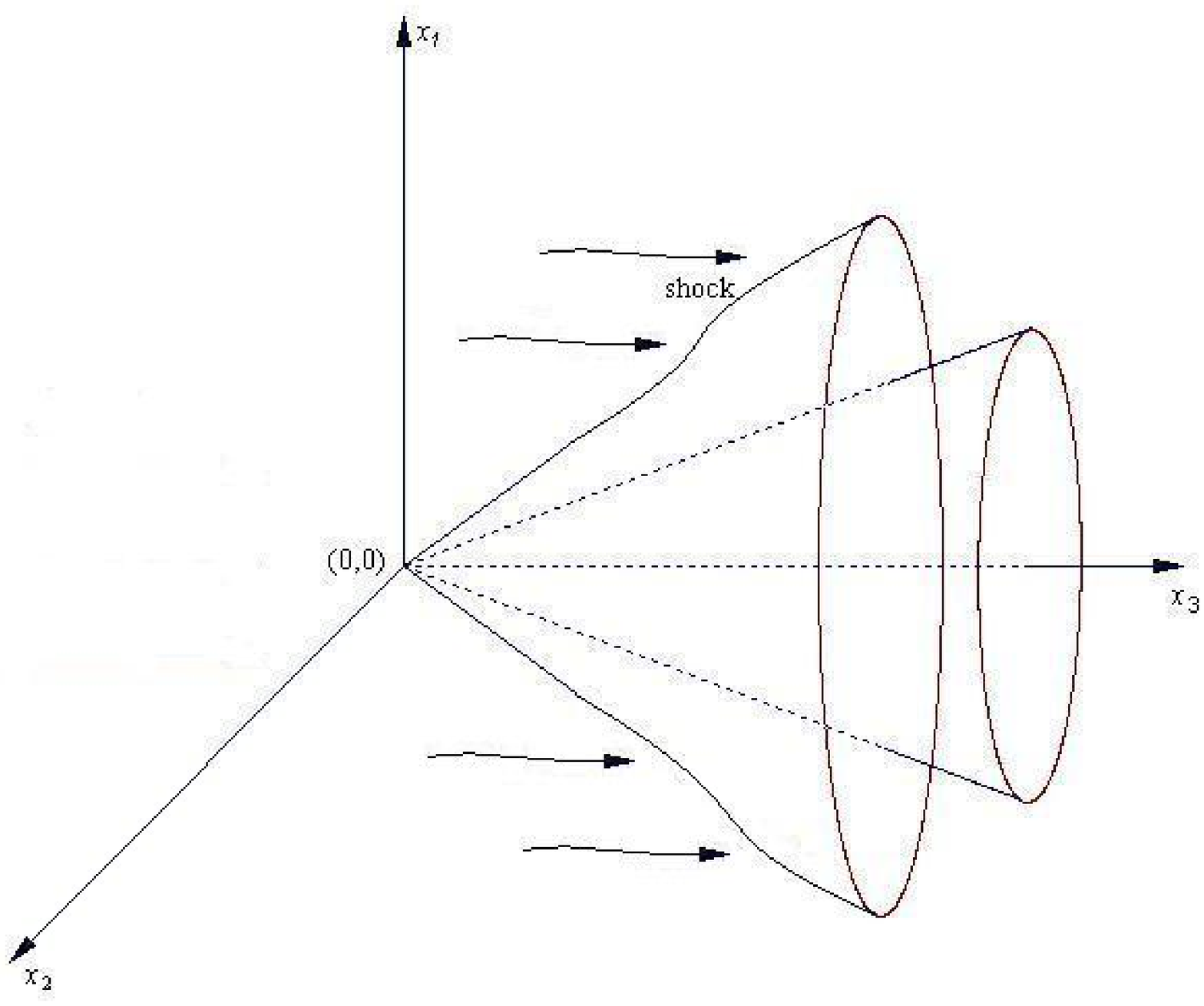}
$$
\centerline{\bf Figure 2. A multi-dimensionally perturbed supersonic
flow past a sharp circular cone}

\vskip 0.6 true cm

In this paper, we will use the potential equation to describe the
motion of a supersonic polytropic gas (this model is also favored in
[22]--[23], [30]), where polytropic gas means that the pressure $P$
and the density~$\rho$ of the gas are related by the equation of state
$P=A\rho^{\g}$, with $A>0$ being a constant and the adiabatic
constant~$\g$ satisfying $1<\g<3$ (for air, $\g\simeq 1.4$). Let
$\Phi(x)$ be a potential of the velocity $u=(u_1,u_2,u_3)$, i.e.,
$u_i=\p_i\Phi$. Then it follows from Bernoulli's law that
$$\f{1}{2}|\na\Phi|^2+h(\rho)=C_0.\tag1.1$$ Here,
$h(\rho)=\ds\f{c^2(\rho)}{\g-1}$ is the specific enthalpy,
$c(\rho)=\sqrt{P'(\rho)}$ is sound speed, $\na=(\p_1,\p_2,\p_3)$, and
$C_0=\ds\f{1}{2}q_0^2+h(\rho_0)$ is Bernoulli's constant which is
determined by the incoming uniform supersonic flow at minus infinity
with velocity $(0,0,q_0)$ and density $\rho_0>0$.

By (1.1) and the implicit function theorem, in view of
$h'(\rho)=\ds\f{c^2(\rho)}{\rho}>0$ for $\rho>0$, the density
$\rho(x)$ can be expressed as
$$\rho=h^{-1}\biggl(C_0-\f{1}{2}|\na\Phi|^2\biggr)
\equiv H(\na\Phi).\tag1.2$$
Substituting (1.2) into the equation $\dsize\sum_{i=1}^3\p_i(\rho
u_i)=0$, which expresses the conservation of mass, yields
$$\dsize\sum_{i=1}^3\p_i\bigl(H(\na\Phi)\p_i\Phi\bigr)=0.\tag1.3$$
More intuitively, for any $C^2$-solution $\Phi$, (1.3) can be
rewritten as a second-order quasilinear equation,
$$
\ds\sum_{i=1}^3((\p_i\Phi)^2-c^2)\p_i^2\Phi+2 \ds\sum_{1\le i<j\le
3}\p_i\Phi\p_j\Phi\p_{ij}^2\Phi=0;\tag 1.4
$$
here $c=c(\rho)=c(H(\na\Phi))$. Note that (1.4) is strictly
hyperbolic with respect to the $x_3$-direction in case $u_3>c(\rho)$
holds.

For the geometry of the conic surface, it is convenient to work in the
cylindrical coordinates $(z, r,\th)$, where
$$x_1=r\cos\th,\quad x_2=r\sin\th,\quad x_3=z,\tag1.5$$
$r=\sqrt{x_1^2+x_2^2}$, and $0\le\th\le 2\pi$.
Under the change of coordinates (1.5), Eq.~(1.4) becomes
$$
\align &\bigl((\p_z\Phi)^2-c^2\bigr)\p_z^2\Phi+
\bigl((\p_r\Phi)^2-c^2\bigr)\p_r^2\Phi+\f{1}{r^2}\biggl(\f{(\p_{\th}\Phi)^2}{
 r^2}-c^2\biggr)\p_{\th}^2\Phi
+2\p_z\Phi\biggl(\p_r\Phi\p_{zr}^2\Phi+\f{1}{r^2}\p_{\th}
\Phi\p_{z\th}^2\Phi\biggr)\\
&\qquad +\f{2}{r^2}\p_r\Phi\p_{\th}\Phi\p_{r\th}^2\Phi-\f{1}{r}\p_r\Phi
\biggl(\f{(\p_{\th}\Phi)^2}{r^2}+c^2\biggr)=0.\tag1.6
\endalign
$$

Let $\Phi^-(z,r,\th)$ and $\Phi^+(z,r,\th)$ denote the velocity
potential for the flow ahead and past the resulting shock front
$r=\chi(z,\th)$, respectively, where $\chi(0,\th)=0$.  Then (1.6)
splits into two equations. That is, $\Phi^{\pm}(z,r,\th)$ satisfy
the following equations in their corresponding domains,
$$
\align &\bigl((\p_z\Phi^-)^2-(c^-)^2\bigr)\p_z^2\Phi^- +
\bigl((\p_r\Phi^-)^2-(c^-)^2\bigr)\p_r^2\Phi^-
+\f{1}{r^2}\biggl(\f{(\p_{\th}\Phi^-)^2}{r^2}-(c^-)^2\biggr)\p_{\th}^2\Phi^-
\\
&\qquad+2\p_z\Phi^-\biggl(\p_r\Phi^-\p_{zr}^2\Phi^-
+\f{1}{r^2}\p_{\th}\Phi^-\p_{z\th}^2\Phi^-\biggr)+\f{2}{r^2}\p_r\Phi^-\p_{\th}
\Phi^-\p_{r\th}^2\Phi^-\\
&\qquad-\f{1}{r}\p_r\Phi^-\biggl(\f{(\p_{\th}\Phi^-)^2}{r^2}
+(c^-)^2\biggr)=0\quad\text{in $\O_-$} \tag1.7
\endalign
$$
and
$$
\align &\bigl((\p_z\Phi^+)^2-(c^+)^2\bigr)\p_z^2\Phi^+ +
\bigl((\p_r\Phi^+)^2-(c^+)^2\bigr)\p_r^2\Phi^+
+\f{1}{r^2}\biggl(\f{(\p_{\th}\Phi^+)^2}{r^2}-(c^+)^2\biggr)\p_{\th}^2\Phi^+
\\
&\qquad+2\p_z\Phi^+\biggl(\p_r\Phi^+\p_{zr}^2\Phi^+
+\f{1}{r^2}\p_{\th}\Phi^+\p_{z\th}^2\Phi^+\biggr)+\f{2}{r^2}
\p_r\Phi^+\p_{\th}\Phi^+\p_{r\th}^2\Phi^+\\
&\qquad-\f{1}{r}\p_r\Phi^+\biggl(\f{(\p_{\th}\Phi^+)^2}{r^2}
+(c^+)^2\biggr)=0\quad\text{in $\O_+$;} \tag1.8
\endalign
$$
here $c^{\pm}=c\bigl(H(\na\Phi^{\pm})\bigr)$, $\O_-=\{(z,r,\th)\colon
r>\chi(z,\th), 0\le\th\le 2\pi, z>0\}$, and $\O_+=\{(z,r,\th)\colon b_0
z\le r<\chi(z,\th), 0\le\th\le 2\pi, z>0\}$.

On the conic surface $r=b_0 z$, $\Phi^+$ satisfies the boundary
condition
$$
\p_r\Phi^+-b_0\p_z\Phi^+=0 \quad \text{on $r=b_0 z$,}\tag1.9
$$
while on the conic shock $\Gamma=\{ r=\chi(z,\th)\}$, by Eq.~(1.3) and
the change of coordinates (1.5), the Rankine-Hugoniot condition
becomes
$$
[H(\na\Phi)\p_r\Phi]-[H(\na\Phi)\p_z\Phi]\p_z\chi=\f{1}{r^2}
[H(\na\Phi)\p_{\th}\Phi]\p_{\th}\chi
\quad \text{on $\Gamma$}.\tag1.10$$

Moreover, the potential $\Phi(z,r,\th)$ is continuous across the
shock, i.e.,
$$\Phi^+\bigl(z,\chi(z,\th),\th\bigr)=\Phi^-\bigl(z,\chi(z,\th),\th\bigr)
\quad \text{on $\Gamma$}.\tag1.11$$

Furthermore, we impose initial conditions on $\Phi^-(z,r,\th)$,
$$\Phi^-(0,r,\th)
=\ve \Phi_0^-(r,\th),\qquad \p_z\Phi^-(0,r,\th) =q_0+\ve
\Phi_1^-(r,\th),\tag1.12$$
where $\ve>0$ is a small parameter, $q_0>c(\rho_0)$, and
$\Phi_0^-(r,\th), \Phi_1^-(r,\th)\in C_0^\infty\bigl((0,l)\times [
0,2\pi]\bigr)$ for some fixed number $l>0$.

\medskip

The main result states:

\smallskip

{\bf Theorem 1.1.} {\it For small $b_0>0$ and a sufficiently large
speed $q_0$, there exists a small constant $\ve_0>0$
depending on $q_0$, $\rho_0$, $b_0$, and $\g$ such that problem
\rom{(1.7)-(1.8)} together with \rom{(1.9)-(1.12)} possesses a
global $C^{\infty}$ supersonic shock solution
$(\Phi^{\pm}(z,r,\th),\chi(z,\th))$ for any $\ve<\ve_0$. Moreover,
$(\na_x\Phi^+,\ds\f{\chi(z,\th)}{z})$ approaches the corresponding
quantities for the incoming uniform supersonic flow $(0,0,q_0)$ with
density $\rho_0$ past the sharp circular cone $r=b_0z$ with rate
$(1+z)^{-m_0}$ for any positive number $m_0<\f12$.}

\medskip

{\bf Remark 1.1.} {\it The various smallness assumptions
in\/ \rom{Theorem 1.1} can be expressed by saying that
$$0<\ve \ll \min\biggl\{\f{1}{b_0^2}\c, \f{1}{b_0^2}\b\biggr\}
\ll b_0^2\quad \text{\rom{and}}\quad b_0\ll b^*,
$$
where $b_*$ is the critical value given in\/ \rom{Remark 2.1} below.}

\medskip

{\bf Remark 1.2.} {\it As in\/ \rom{[5]}, \rom{[9-10]},
and \rom{[31]}, where suitable symmetry assumptions on the incoming
supersonic flow with small perturbation or artificial boundary
condition on the conic surface were imposed, we emphasize that also in
the case treated in this paper there are no discontinuities for the
weak solution $\Phi(x)=\Phi^+(x)$ in $\O_+$ and $\Phi^-(x)$ in $\O_-$,
but the main multi-dimensional conic shock front. This means that the
supersonic conic shock is structurally stable in the whole space for a
polytropic gas and arbitrary perturbations. This agrees with
observations from physical experiment and numerical simulations.}

\medskip

{\bf Remark 1.3.} {\it The nonlinear hyperbolic equation\/ \rom{(1.4)}
is actually a second-order quasilinear wave equation in two space
dimensions if one regards $x_3$ as time, as the flow is supersonic in
$x_3$-direction. By a direct verification, one sees that\/ \rom{(1.4)}
does not fulfill the ``null-condition'' put forward in\/ \rom{[6]}
and \rom{[19]}. Therefore, in terms of the results
of \rom{[1]}, \rom{[12--13]}, \rom{[15--17]}, \rom{[25--26]},
and \rom{[28]}, if there was no shock for the solution
to \rom{Eq.~(1.4)}, then the classical solution to
\rom{(1.4)} would blow up in finite time. Thus, the result of\/
\rom{Theorem~1.1} asserts that the multi-dimensional supersonic
shock absorbs all possible compressions in the flow and prevents the
flow from forming further shocks as well as other singularities.}

\medskip

{\bf Remark 1.4.} {\it As energy estimates fail to hold
in\/ \rom{BV}-spaces for multi-dimensional hyperbolic equations and
systems as shown in\/ \rom{[27]}, the method used
in\/ \rom{[21]} \rom(which was the Glimm scheme for a spherically
symmetrically perturbed conic surface\/\rom) cannot be applied to the
genuinely multi-dimensional problem treated here.}

\medskip

{\bf Remark 1.5.} {\it It was indicated in \rom{[7, pages 317--318 and
414]} that if a supersonic steady flow approaches from minus infinity
and hits a sharp cone in direction of its axis, then it follows from
the Rankine-Hugoniot conditions and the physical entropy condition, by
an application of the method of the apple
curve\/ \rom(see\/ \rom{Figure~3} below\/\rom) that there possibly
occur a weak shock and a strong shock attached to the tip of the
cone. These shocks are supersonic and transonic, respectively. It was
frequently stated that the strong shock is unstable and that,
therefore, only the weak shock is present in real
situations. In \rom{[32-34]}, the global instability of an attached
strong conic shock in the whole space was systematically
studied\/ \rom(in this case, the corresponding subsonic potential
equation is nonlinear elliptic and the steady Euler system becomes
elliptic-hyperbolic\/\rom) which especially showed that a global
strong conic shock is actually unstable as long as the perturbation of
the sharp circular cone satisfies some suitable assumptions. On the
other hand, from the result in this paper, one infers the global
stability of a supersonic conic shock. Consequently, in regard to
the global stability or instability of weak and strong conic shocks,
these results give a partial answer to the problem of Courant and
Friedrichs\/ \rom{[7]} stated above.}
$$
\epsfysize=76mm \epsfbox{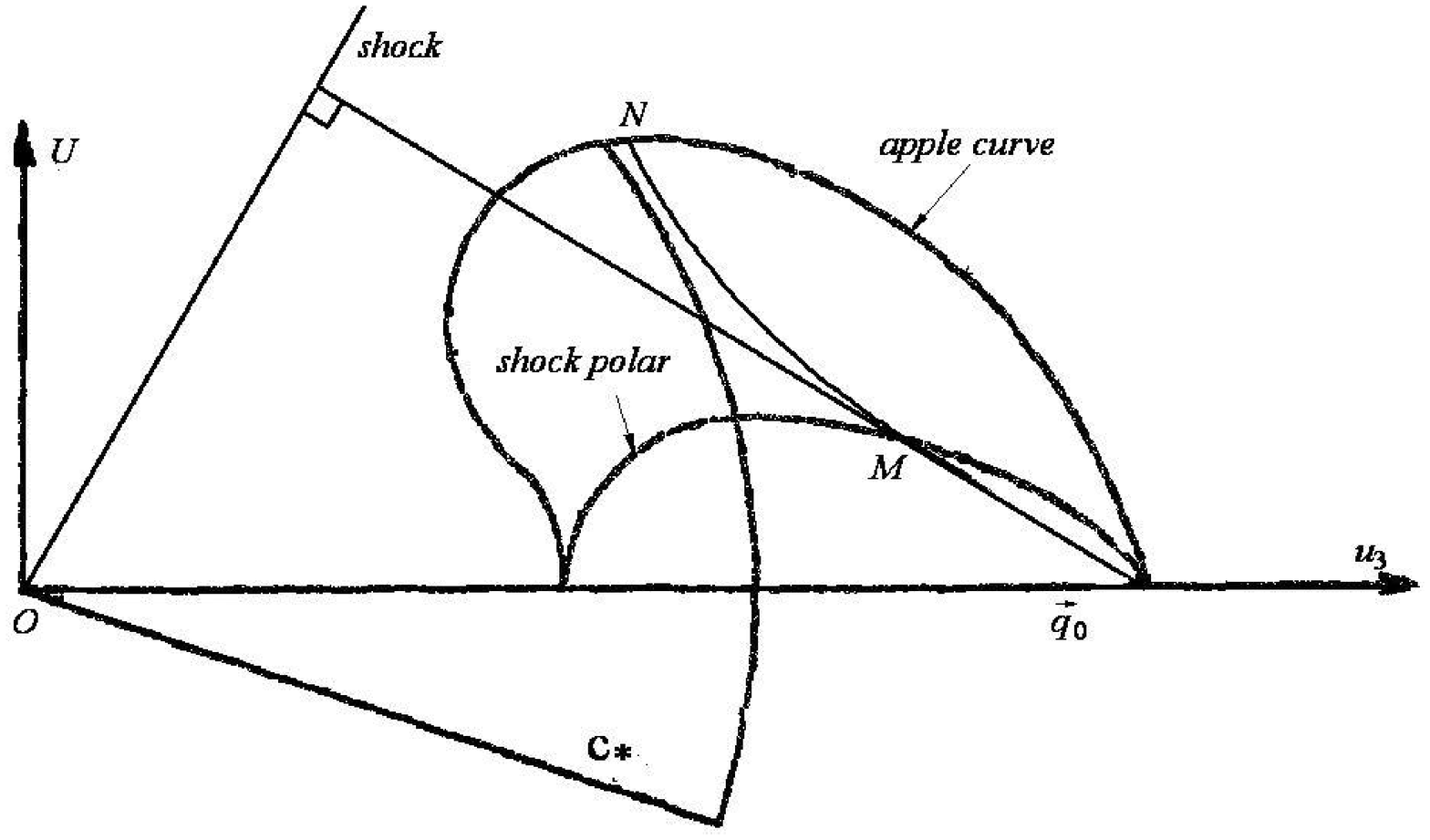}
$$
\centerline{\bf Figure 3. Apple curve showing all possible end
velocity of a conic shock.}

\vskip 0.4 true cm

{\bf Remark 1.6.} {\it For large $q_0$, the incoming flow is called
hypersonic. The famous independence principle for large Mach
numbers\/ \rom(that there exists a stable limit state for a hypersonic
flow as the Mach number goes to infinity\/\rom) is likewise
illustrated by\/ \rom{Theorem~1.1} for a hypersonic gas past a sharp
cone. For more physical properties of hypersonic flow,
see\/ \rom{[2]}, \rom{[8]}, and \rom{[29]}.}

\smallskip

{\bf Remark 1.7.} {\it In case the conic surface is also perturbed,
where the multi-dimensional perturbation is small and of compact
support\/ \rom(possibly including a compact perturbation near the
vertex of the cone, as local existence of such a conic shock was
established in\/ \rom{[4])} or decays sufficiently fast when $z$ goes
to infinity, a result analogous to\/ \rom{Theorem 1.1} remains in
force, where the proofs of this paper still work.}

\smallskip

Let us mention some work which is directly related to this paper. In
[5] and [9], under the assumptions of an incoming uniform supersonic
flow and that the angle of the spherically symmetrically curved conic
body is sufficiently small (smaller than the critical angle which
guarantees that the supersonic shock is attached), it was shown that a
spherically symmetric supersonic conic shock exists globally past the
conic body for a supersonic polytropic flow.  Z.~Xin and H.~Yin
established in [31] the global existence of a multi-dimensional
supersonic conic shock for an incoming uniform supersonic flow past a
generally curved sharp cone under an artificial boundary condition --
the Dirichlet condition for the potential on the conic surface
(physically, this kind of boundary condition means that the conic body
is perforated or porous). It should be emphasized here that the
Dirichlet boundary condition for the potential in [31] played a
crucial role in deriving {\it a priori\/} energy estimates and further
obtaining the global existence. It means that the Poincar\'e
inequalities are available on the shock surface and in the interior of
the downstream domain. By applying the Glimm scheme, in the case of a
spherically symmetrically curved cone, W.-C.~Lien and T.-P.~Liu in
[21] obtained the global existence of a weak solution and the
long-distance asymptotic behavior under suitable restrictions on the
large Mach number, the sharp vertex angle, and the shock strength. The
main interest here is to establish the global existence of a genuinely
multi-dimensional supersonic conic shock for a perturbed supersonic
polytropic gas past an infinitely long cone with a sharp angle when
the speed of the incoming flow is large. Especially, we remove the key
assumption of spherically symmetry on the perturbed supersonic flow
which was assumed in [5], [9], and [21] and which was essential in the
proofs there.

Let us also comment on the proof of Theorem 1.1.  In order to prove
Theorem 1.1, we intend to use continuous induction to establish {\it a
priori\/} estimates of the solution and its derivatives.  To achieve
this objective, as in [5], [11], and [15], we need to derive global
weighted energy estimates for the linearized problem (1.7)--(1.8) with
(1.9)--(1.12). Based on such estimates, one then obtains the global
existence, stability, and the asymptotic behavior of the shock
solution to the perturbed nonlinear problem. The key ingredients in
the analysis to obtain weighted energy estimates are an appropriate
multiplier and a new Hardy-type inequality on the shock
surface. Finding a suitable multiplier is a hard task for the
following reasons: First, to obtain the global existence requires to
establish global estimates, independent of $z$, of the potential
function and its derivatives on the boundaries as well as in the
interior of the downstream domain. This implies strict constraints on
the multiplier and makes the computations delicate and
involved. Secondly, as our background solution is self-similar in a
downstream domain and strongly depends on the location of the boundary
of the cone, the angle at the vertex of the cone, the Mach number of
the incoming flow, and the equation of state of the gas under
consideration, one needs to take some measures to simplify the
coefficients of the nonlinear equation together with the corresponding
nonlinear boundary conditions so that the procedure to find the
multiplier becomes managable. Thirdly, for the multi-dimensional case,
the Neumann-type boundary condition (1.9) fulfilled by $\Phi^+$
introduces additional difficulties compared to [5] and [31], where [5]
only treats the case of a spherically symmetric conic shock with
Neumann-type boundary condition on the conic surface, while [31]
treats the case of an artificial Dirichlet-type boundary condition for
the potential on a multi-dimensionally perturbed conic surface. The
latter plays a key role in the analysis of [31], as the corresponding
Poincar\'e inequality is available on the shock surface and the
interior of the downstream domain, respectively, while this is not the
case in the problem treated here. Thanks to some delicate analysis
accompanied by a new Hardy-type inequality derived by making full use
of the special structure of the shock boundary conditions (i.e., the
sizes as well as the signs of the coefficients in (3.6)--(3.7)), we
finally overcome all these difficulties and obtain a uniform estimate
of $\|\na_x\Phi^+\|_{L^2(\O_+)}$. From this, higher-order estimates of
$\na_x\Phi^+$ can be established by using Klainerman vector field and
commutator arguments together with a careful verification that
suitably chosen higher-order derivatives of the solution satisfy the
Neumann-type boundary condition on the conic surface.  This eventually
establishes Theorem~1.1.

The paper is organized as follows: In \S2, we derive some basic
estimates of the background self-similar solution in case of an
incoming hypersonic flow, which are required to treat the
linearization of the nonlinear problem and for the construction of the
multiplier. In \S3, we reformulate problem (1.7)--(1.12) by
decomposing its solution as a sum of the background solution and a
small perturbation $\dot\vp$ so that its linearization can be studied
in a convenient way. In \S4, we first establish a uniform weighted
energy estimate for the corresponding linear problem, where also an
appropriate multiplier is constructed. Based on such an energy
estimate, we obtain a uniform weighted energy estimate of
$\na_x\dot\vp$ for the nonlinear problem through establishing a new
Hardy-type inequality. In \S5, by the estimates derived in $\S 4$, we
continue to establish uniform higher-order weighted estimates of
$\na_x\dot\vp$. In \S6, the proof of Theorem~1.1 is eventually
completed by utilizing Sobolev's embedding theorem and continuous
induction. Some lengthy computations are carried out in
an appendix.

\smallskip

In what follows, we will use the following conventions:

$\centerdot$ $C$ stands for a generic positive constant which does not
depend on any quantity except the adiabatic constant $\gamma$
($1<\g<3$).

$\centerdot$ $C(\cdot)$ represents a generic positive constant
which depends on its argument (or arguments).

$\centerdot$ $O(\cdot)$ means that $|O(\cdot)|\leq C|\cdot|$ holds
true. In particular, $O(\ve)$ abbreviates $|O(\ve)|\leq
C(b_0,q_0)\ve$.

$\centerdot$ $dS$ stands for the surface measure in the
corresponding surface integral.

\vskip 0.4 true cm \centerline{\bf \S2. Analysis of the
self-similar background solution} \vskip 0.4 true cm

In this section, we will provide, with more details than in [5], 
properties of the background solution when the Mach number of
the incoming supersonic flow is large and the supersonic shock is
attached.  These properties will be applied again and again in the
later analysis of \S3--\S5.

Following the illustrations of [7, page 407], the supersonic conic
shock phenomenon for an incoming supersonic flow past a sharp circular
cone is described as follows: Suppose that there is a uniform
supersonic flow $(0,0,q_0)$ with constant density $\rho_0>0$ which
approaches from minus infinity. Let the flow hit the circular cone
$\{(r,z)\colon r\le b_0 z,\, z\ge 0\}$ in the direction of the
$z$-axis. Then there exists a critical value $b^*>0$, which is
determined by the parameters of the incoming flow, such that there
occurs a supersonic conic shock $r=s_0 z$, $s_0>b_0$, attached to the
tip of the cone whenever $b_0<b^*$ holds true. Moreover, the solution
to (1.3) with (1.1) past the shock surface is self-similar, that is,
in cylindrical coordinates $(z, r,\th)$, the density and velocity
between the shock front and the conic surface are of the form
$\rho=\rho(s)$, $u_1=u_r(s)\ds\frac{x_1}{r}$,
$u_2=u_r(s)\ds\frac{x_2}{r}$, and $u_3=u_z(s)$, where
$s=\ds\frac{r}{z}$. In this case, Eq.~(1.3) with (1.1) can be reduced
to the following nonlinear ordinary differential system:
$$
\cases
&\ds\rho'(s)=-\frac{\rho u_r(su_z-u_r)}{s\bigl((1+s^2)
c^2(\rho)-(su_z-u_r)^2\bigr)},\\
&\ds u_r'(s)=-\frac{c^2(\rho) u_r}{s\bigl((1+s^2)
c^2(\rho)-(su_z-u_r)^2\bigr)},\\
&\ds u_z'(s)=\frac{c^2(\rho) u_r}{(1+s^2)c^2(\rho)-(su_z-u_r)^2},
\endcases\quad
\text{for $b_0\leq s\leq s_0$.}\tag2.1
$$
As explained in [5] or [7], for the denominator it holds
$(1+s^2)c^2(\rho)-(su_z-u_r)^2>0$ for $b_0\leq s\leq s_0$, which
implies that system (2.1) makes sense.

On the shock front $r=s_0 z$, it follows from the Rankine-Hugoniot
conditions and Lax's geometric entropy conditions on the $2$-shock
that
$$
\cases &[\rho u_r]-s_0[\rho u_z]=0,\\
&[u_z]+s_0[u_r]=0
\endcases\tag2.2
$$
and
$$
\cases
&\lambda_1(s_0)<s_0<\lambda_2(s_0),\\
&\ds\frac{c(\rho_0)}{\sqrt{q_0^2-c^2(\rho_0)}}<s_0,
\endcases\tag2.3
$$
where $$\ds\lambda_{1,2}(s)=\frac{u_r(s)u_z(s)\mp
c(\rho(s))\sqrt{u_r^2(s)+u_z^2(s)-c^2(\rho(s))}}{u_z^2(s)-c^2(\rho(s))}.\tag
2.4$$
Additionally, the flow satisfies the condition
$$u_r(b_0)=b_0 u_z(b_0)\tag2.5$$
on the fixed boundary $s=b_0$.

As indicated in [7, pages 411--414] or [18], the boundary value
problem (2.1)--(2.4) can be solved by the shooting method as well as
by the method of the apple curve. More specifically, for any given
$b_0>0$, which is smaller than the critical value $b^*$, one can
determine the solution to (2.1)--(2.4) by finding the integral curve
of $\ds\f{du_r}{du_z}=-\ds\f{u_z}{u_r}$ from the intersection point of
the apple curve with the ray $u_r=b_0u_z$ to some point at the shock
polar (see Figure~3 above or [7, Fig.~8 on page~414]). In this paper,
such a supersonic solution past the shock is called the background
solution.

\medskip

For large $q_0$, some detailed properties of the background
solution are as follows:

\smallskip

{\bf Lemma 2.1.} {\it For $q_0$ large enough, $1<\gamma<3$, and
$0<b_0<b_*=\sqrt{\ds\f{1}{2}(\sqrt{\f{\g+7}{\g-1}}-1)}$, one has, for
$b_0\le s\le s_0$,

\rom{(i)} $s_0=b_0\biggl(1+ O(\b)\biggr)$.

\rom{(ii)} $0\le su_z(s)-u_r(s)\le O(\a)$.

\rom{(iii)} $u_r(s)=\ds\f{b_0q_0}{1+b_0^2}\biggl(1+O(\b)\biggr)$.

\rom{(iv)} $u_z(s)=\ds\f{q_0}{1+b_0^2}\biggl(1+ O(\b)\biggr)$.

\rom{(v)}
$\rho(s)=\biggl(\ds\f{\g-1}{2A\g(1+b_0^2)}\biggr)^{\f{1}{\g-1}} (b_0
q_0)^{\f{2}{\g-1}}\biggl(1+O(\c)+O(\b)\biggr)$.

\rom{(vi)}
$q^2(s)-c^2(\rho(s))=q_0^2\biggl(\ds\frac{2-(\gamma-1)b_0^2}{2(1+b_0^2)}\biggr)\biggl((1+
O(\c) + O(\b)\biggr)$,

\qquad\qquad\qquad\qquad\qquad\qquad\qquad\qquad here and below
$q^2(s)=u_r^2(s)+u_z^2(s)$.

\rom{(vii)}
$u_z^2(s)-c^2(\rho(s))=\ds\f{1-\f{\g-1}{2}b_0^2(1+b_0^2)}{(1+b_0^2)^2}q_0^2\biggl(1+
O(\c)+ O(\b)\biggr)>0$.

\rom{(viii)}
$(1+s^2)c^2(\rho(s))-(su_z(s)-u_r(s))^2=\ds\f{\g-1}{2\g}(b_0q_0)^2\biggl(1+O(\c)+O(\b)\biggr)>0$.
}\vskip 0.2 true cm

\smallskip

{\bf Remark 2.1.} {\it
$b_*=\sqrt{\ds\f{1}{2}(\sqrt{\f{\g+7}{\g-1}}-1)}$ is a root of the
quartic equation $1-\f{\g-1}{2}b_0^2(1+b_0^2)=0$. This guarantees
that the flow across the shock is still supersonic in direction of $z$
for $b_0<b_*$ and large $q_0$. It can be derived from the
expression in\/ \rom{Lemma 2.1 (vii)}.}

\vskip 0.2 true cm

{\bf Proof.} Set $\rho^{+}=\lim\limits_{s\to s_{0}-}\rho(s)$,
$u_r^{+}=\lim\limits_{s\to s_{0}-}u_r(s)$, $u_{z}^{+}=\lim\limits_{s\to
s_{0}-}u_{z}(s)$, and $\al=\ds\f{\rho^+}{\rho_0}$.
It follows the Rankine-Hugoniot conditions (2.2) that
$$\cases
&\ds u_z^{+}=\f{q_0}{1+s_0^2}\biggl(1+\f{s_0^2}{\al}\biggr),\\
&\ds u_r^{+}=\f{s_0 q_0}{1+s_0^2}\biggl(1-\f{1}{\al}\biggr).
\endcases\tag 2.6$$
Substituting (2.6) into Bernoulli's law (1.1) yields
$$\f{A\g}{\g-1}\biggl((\rho^{+})^{\g+1}-\rho_0^{\g-1}(\rho^{+})^2\biggr)
+\f{s_0^2
q_0^2}{2(1+s_0^2)}\biggl(\rho_0^2-(\rho^{+})^2\biggr)=0.\tag 2.7$$

Set
$$F(x)=\f{A\g}{\g-1}x^{\g+1}-\f{A\g}{\g-1}\rho_0^{\g-1}x^2+\f{s_0^2
q_0^2}{2(1+s_0^2)}(\rho_0^2-x^2).$$
Then (2.7) implies that $F(\rho_0)=F(\rho^+)=0$.
It follows from a direct computation that, for $x\in (0,\rho_0)$,
$$\align
F'(x)&=\f{A\g(\g+1)}{\g-1}x^{\g}-\f{2A\g}{\g-1}\rho_0^{\g-1}x-\f{s_0^2
q_0^2}{1+s_0^2}x
=x\biggl(\f{A\g(\g+1)}{\g-1}x^{\g-1}-\f{2A\g}{\g-1}\rho_0^{\g-1}-\f{s_0^2
q_0^2}{1+s_0^2}\biggr)\\
&\leq
x\biggl(\f{A\g(\g+1)}{\g-1}\rho_0^{\g-1}-\f{2A\g}{\g-1}\rho_0^{\g-1}-\f{s_0^2
q_0^2}{1+s_0^2}\biggr)
=x\biggl(c^2(\rho_0)-\f{s_0^2 q_0^2}{1+s_0^2}\biggr).
\endalign$$
Due to Lax's geometric entropy conditions (2.3), $F'(x)<0$ for $x\in
(0,\rho_0)$. Together with $F(\rho_0)=F(\rho^+)=0$, this yields
$\rho^+>\rho_0$.

In this case, (2.7) is equivalent to
$$\al^2\f{\al^{\g-1}-1}{\al^2-1}=\f{(\g-1)s_0^2
q_0^2}{2A\g(1+s_0^2)\rho_0^{\g-1}},\tag 2.8$$ where $\al>1$.  Since
the left hand side of (2.8) is bounded if $\al>1$ is bounded, for
$q_0$ is large, $\al$ is also large. From this fact, one has
$$\al=\f{1}{\rho_0}\biggl(\f{\g-1}{2A\g(1+s_0^2)}
\biggr)^{\f{1}{\g-1}}(b_0q_0)^{\f{2}{\g-1}}\biggl(1+O(\c)+O(\b)\biggr).$$
Substituting this into (2.6) yields
$$\cases
&\ds u_{z}^{+}=\f{q_0}{1+s_0^2}\biggl(1+O(\b)\biggr),\\
&\ds u_{r}^{+}=\f{s_0 q_0}{1+s_0^2}\biggl(1+O(\b)\biggr).
\endcases\tag 2.9$$

Moreover, from $u_r'(s)<0$ and $u_z'(s)>0$ for $s\in [b_0,s_0]$, one has
$$u_r^{+}\leq u_r(s)\leq u_r(b_0)=b_0 u_{z}(b_0)\leq b_0 u_{z}(s)\leq b_0
u_{z}^+.\tag 2.10$$
Combining (2.9) with (2.10) yields (i) and further (iii)--(iv).

(ii) comes from (2.9) and the fact that
$$0=b_0 u_z(b_0)-u_r(b_0)\leq s u_z(s)-u_r(s)\leq s_0
u_z^+-u_r^+\quad\text{for $b_0\leq s\leq s_0$};$$
here $(s u_z(s)-u_r(s))'\ge 0$ for $b_0\leq s\leq s_0$ has been applied.

Furthermore, Bernoulli's law (1.1) and (i)--(iv) show by a direct
computation that (v)--(viii) hold. Therefore, the proof of Lemma 2.1
is complete. \qed

\medskip

{\bf Lemma 2.2.} {\it Under the assumptions of\/ \rom{Lemma 2.1}, one
has, for $b_0\le s\le s_0$,

\rom{(i)} $\la_1(s)-s=\ds\f{\sqrt{\g-1}(1+b_0^2)b_0\bigl(\sqrt{\g-1}
b_0^2-\sqrt{2-(\g-1)b_0^2}\bigr)}{2-(\g-1)b_0^2(1+b_0^2)}
\bigl(1+O(\c)+O(\b)\bigr)<0$.

\rom{(ii)} $\la_2(s)-s=\ds\f{\sqrt{\g-1}(1+b_0^2)b_0\bigl(\sqrt{\g-1}
b_0^2+\sqrt{2-(\g-1)b_0^2}\bigr)}{2-(\g-1)b_0^2(1+b_0^2)}
\bigl(1+O(\c)+O(\b)\bigr)>0$.

\rom{(iii)}
$u_r'(s)=-\ds\frac{q_0}{(1+b_0^2)^2}\biggl(1+O(\c)+O(\b)\biggr)$.

\rom{(iv)}
$u_z'(s)=\ds\frac{b_0q_0}{(1+b_0^2)^2}\biggl(1+O(\c)+O(\b)\biggr)$.

\rom{(v)} $|\rho'(s)|\le O(\f{1}{b_0})$.}

{\bf Proof.} (i)--(v) can be directly derived from (2.1), (2.4), and
Lemma 2.1. For the reader's convenience, we provide the detailed
computation for (i) as an example.

It follows from (2.4), (iii)--(v) in Lemma 2.1, and a direct
computation that
$$\align
\lambda_1(s)&=\ds\f{u_{r}(s)u_{z}(s)-c(\rho(s))\sqrt{u_r^2(s)+u_z^2(s)
-c^2(\rho(s))}}{u_z^2(s)-c^2(\rho(s))}\\
&=\ds\f{\f{b_0 q_0}{1+b_0^2}\f{q_0}{1+b_0^2}-b_0
q_0\sqrt{\f{\g-1}{2(1+b_0^2)}}\sqrt{\f{(b_0
q_0)^2}{(1+b_0^2)^2}+\f{q_0^2}{(1+b_0^2)^2}-\f{\g-1}{2(1+b_0^2)}(b_0
q_0)^2}}{\f{q_0^2}{(1+b_0^2)^2}-\f{\g-1}{2(1+b_0^2)}(b_0
q_0)^2}\biggl(1+O(\c)+O(\b)\biggr)\\
&=\f{2b_0-(1+b_0^2)\sqrt{\g-1}b_0\sqrt{2-(\g-1)b_0^2}}{2-(\g-1)b_0^2(1+b_0^2)}
\biggl(1+O(\c)+O(\b)\biggr).\\
\endalign$$
Due to $s_0=b_0\biggl(1+O(\b)\biggr)$ by (i) in Lemma 2.1,
$\lambda_1(s)-s$ satisfies (i). \qed

\smallskip

{\bf Remark 2.2.} {\it Since the denominator of system\/ \rom{(2.1)}
is positive in the interval $[b_0, s_0]$ \rom(see\/ \rom{Lemma 2.1
(viii))}, one can extend the background solution $(\rho(s),
u_z(s),u_r(s))$ to\/ \rom{(2.1)--(2.3)} and\/ \rom{(2.5)} to the
interval $[b_0, s_0+\tau_0]$ for some small positive constant $\tau_0$
satisfying $0<\tau_0\leq q_0^{-\frac{4}{\gamma-1}}(s_0-b_0)$. In the
following sections, we will denote this extension of the background
solution to $\{(z, r)\colon z>0,\, b_0 z\leq r\leq (s_0+\tau_0)z\}$ by
$(\hat\rho(s), \hat u_z(s), \hat u_r(s))$, where $s=\ds\f{r}{z}$. The
corresponding extension of the potential will be denoted by
$\hat\Phi(s)$.}

\vskip 0.4 true cm \centerline{\bf \S3. Reformulation of
the nonlinear problem} \vskip 0.4 true cm

In this section, we reformulate problem (1.8)--(1.11) by decomposing
its solution as a sum of the background solution and a small
perturbation. Moreover, based on the analysis of the background
solution in Lemmas 2.1 and 2.2, we establish estimates of the
coefficients which appear in the reformulated problem when $q_0$ is
large.

\medskip

We now provide a global existence result for the solution to Eq.~(1.7)
with initial data (1.12) ahead of the shock.

\smallskip

{\bf Lemma 3.1.} {\it \rom{Eq.~(1.7)} with \rom{(1.12)} possesses a
$C^{\infty}$-solution $\Phi^-(z,r,\th)$ in $\O_-=\{(z,r,\th)\colon
r>\chi(z,\th),\, 0\le\th\le 2\pi, z>0\}$.  Moreover,
$\Phi^-(z,r,\th)-q_0z\in C_0^{\infty}(\O_-)$ and, for any $k\in\Bbb
N$, there exists a positive constant $C_k$ independent of $\ve$ such
that
$$\|\Phi^-(z,r,\th)-q_0z\|_{C^k(\O_-)}\leq C_k\ve.
$$}

{\bf Proof.} The quasilinear equation (1.7) is strictly hyperbolic with
respect to the direction of $z$, for the supersonic flow condition
$u_3^->c(\rho^-)$. The initial condition (1.12) means a small
perturbation with compact support away from the origin. Thus, for
large $q_0$, by finite propagation speed and standard Picard iteration
(or see [16]), one easily derives the validity of Lemma 3.1. \qed

\smallskip

Note that there exists a constant $T_0>0$ such that $\Phi^-(z,r,\th)=
q_0 z$ for $z>T_0$. Without loss of generality, below we will assume
that $T_0=1$.

Next, we reformulate the nonlinear problem (1.8)--(1.11).  For
notational convenience, we will neglect all the superscripts ``$+$" in
(1.8)--(1.11) from now on.

Let $\Phi$ be the solution to (1.8)--(1.11) and $\dot{\varphi}$ be the
perturbation of the background solution, that is,
$\dot{\varphi}=\Phi-\hat{\Phi}$; here $\hat{\Phi}$ is given in Remark
2.2. Then, by a direct computation, (1.7) is reduced to:
$$\align
\Cal
L\dot\vp&=f_1(\f{r}{z},\na_x\dot\vp)\p_{z}^2\dot\vp
+f_2(\f{r}{z},\na_x\dot\vp)\p_{zr}^2\dot\vp
+f_3(\f{r}{z},\na_x\dot\vp)\p_{r}^2\dot\vp
+\f{1}{r^2}f_4(\f{r}{z},\na_x\dot\vp)\p_{\th}^2\dot\vp\\
& \qquad +\f{1}{r}f_5(\f{r}{z},\na_x\dot\vp)\p_{z\th}^2\dot\vp
+\f{1}{r}f_6(\f{r}{z},\na_x\dot\vp)
\p_{r\th}^2\dot\vp+\f{1}{r}f_7(\f{r}{z},\na_x\dot\vp)\quad \text{in
$\O_+$},
\tag 3.1\\
\endalign
$$
where
$$
\left\{ \enspace
\aligned\Cal
L\dot\vp&=\p_z^2\dot\vp+2P_1(s)\p_{zr}^2\dot\vp+P_2(s)\p_r^2\dot\vp
-\ds\f{1}{r^2}P_3(s)\p_{\th}^2\dot\vp
+\ds\f{2}{r}P_4(s)\p_z\dot\vp+\ds\f{2}{r}P_5(s)\p_r\dot\vp,\\
P_1(s)&=\ds\f{\hat u_z(s) \hat u_r(s)}{{\hat u}_z^2(s))-c^2(\hat \rho(s))},\\
P_2(s)&=\ds\f{{\hat u}_r^2(s)-c^2(\hat \rho(s))}{{\hat u}_z^2(s))
-c^2(\hat \rho(s))},\\
P_3(s)&=\ds\f{c^2(\hat \rho(s))}{{\hat u}_z^2(s)-c^2(\hat \rho(s))},\\
P_4(s)&=\ds\f{1}{{\hat u}_z^2(s)-c^2(\hat \rho(s))}
\biggl(-\f{\g+1}{2} s^2
\hat u_z(s)\hat u_z'(s)+\f{\g-1}{2}s \hat u_z(s) \hat u_r'(s)\\
&\qquad+s\hat  u_r(s) \hat u_z'(s)+\ds\f{\g-1}{2}\hat
u_z(s)\hat u_r(s)
\biggr),\\
P_5(s)&=\ds\f{1}{{\hat u}_z^2(s)-c^2(\hat \rho(s))}
\biggl(-\f{\g-1}{2}s^2
\hat u_r(s) \hat u_z'(s)+\ds\f{\g+1}{2}s \hat u_r(s) \hat u_r'(s)\\
&\qquad+s \hat u_z(s)
\hat u_z'(s)+\ds\f{\g-1}{2}{\hat u}_r^2(s)
-\ds\f{1}{2}c^2(\hat \rho(s)) \biggr)
\endaligned
\right.\tag 3.2
$$
and
$$\align
f_1(s,\na_x\dot\vp)&=\f{1}{\hat u_z^2(s)-c^2(\hat
\rho(s))}\biggl\{-2\hat u_z(s)\p_{z}\dot\vp-(\p_{z}
\dot\vp)^2-\f{\g-1}{2}\biggl((2\hat u_z(s)+\p_{z}\dot\vp)\p_{z}\dot\vp\\
&\qquad +(2\hat
u_r(s)+\p_r\dot\vp)\p_r\dot\vp+\f{(\p_{\th}\dot\vp)^2}{r^2}\biggr)
\biggr\},\\
f_2(s,\na_x\dot\vp)&=\f{1}{\hat u_z^2(s)-c^2(\hat
\rho(s))}\biggl\{-2\hat u_z(s)\p_{r}\dot\vp-2\hat u_r(s)
\p_{z}\dot\vp-2\p_{z}\dot\vp\p_{r}\dot\vp\biggr\},\\
f_3(s,\na_x\dot\vp)&=\f{1}{\hat u_z^2(s)-c^2(\hat
\rho(s))}\biggl\{-2\hat u_r(s)\p_{r}\dot\vp-(\p_r\dot\vp)^2
-\f{\g-1}{2}\biggl((2\hat u_z(s)+\p_{z}\dot\vp)\p_{z}\dot\vp\\
&\qquad +(2\hat
u_r(s)+\p_{r}\dot\vp)\p_{r}\dot\vp+\f{(\p_{\th}\dot\vp)^2}{r^2}
\biggr)\biggr\},\\
f_4(s,\na_x\dot\vp)&=-\f{1}{\hat u_z^2(s)-c^2(\hat
\rho(s))}\biggl\{\f{(\p_{\th}\dot\vp)^2}{r^2}+\f{\g-1}{2}
\biggl((2\hat u_z(s)+\p_{z}\dot\vp)\p_{z}\dot\vp\\
&\qquad +(2\hat u_r(s)+\p_r\dot\vp)\p_r\dot\vp+\f{(\p_{\th}\dot\vp)^2}{r^2}
\biggr)\biggr\},\\
f_5(s,\na_x\dot\vp)&=-\f{2}{\hat u_z^2(s)-c^2(\hat
\rho(s))}\biggl\{\f{1}{r}\hat u_z(s)\p_{\th}\dot\vp+
\f{1}{r}\p_{z}\dot\vp\p_{\th}\dot\vp\biggr\},\\
f_6(s,\na_x\dot\vp)&=-\f{2}{\hat u_z^2(s)-c^2(\hat
\rho(s))}\biggl\{\f{1}{r}\hat u_r(s)\p_{\th}\dot\vp
+\f{1}{r}\p_r\dot\vp\p_{\th}\dot\vp\biggr\},\\
\endalign
$$

$$
\align
f_7(s,\na_x\dot\vp)&=\f{1}{\hat u_z^2(s)-c^2(\hat
\rho(s))}\biggl\{s^2 \hat
u_z'(s)\biggl(\f{\g+1}{2}(\p_z\dot\vp)^2+\f{\g-1}{2}(\p_r\dot\vp)^2
+\f{\g-1}{2}\f{(\p_{\th}\dot\vp)^2}{r^2}\biggr)\\
&\qquad -s\hat
u_r'(s)\biggl(\f{\g-1}{2}(\p_z\dot\vp)^2+\f{\g+1}{2}(\p_r\dot\vp)^2
+\f{\g-1}{2}\f{(\p_{\th}\dot\vp)^2}{r^2}\biggr)-2s\hat u_z'(s)\p_z
\dot\vp\p_r\dot\vp\\
&\qquad -\hat
u_r(s)\biggl(\f{\g-1}{2}(\p_z\dot\vp)^2+\f{\g-1}{2}(\p_r\dot\vp)^2
+\f{\g-3}{2}\f{(\p_{\th}\dot\vp)^2}{r^2}\biggr)\\
&\qquad +\p_r\dot\vp\biggl(\f{(\p_{\th}\dot\vp)^2}{r^2}-\f{\g-1}{2}\bigl((2\hat
u_z(s)+\p_z\dot\vp)\p_z\dot\vp+ (2\hat
u_r(s)+\p_r\dot\vp)\p_r\dot\vp+\f{(\p_{\th}\dot\vp)^2}{r^2}\bigr)
\biggr)\biggr\},
\endalign$$
where $s=\ds\f{r}{z}$.

On $r=b_0z$, one has
$$\p_r\dot\vp=b_0\p_z\dot\vp.\tag3.3$$

On the free boundary $r=\chi(z,\th)$, by the continuity condition
(1.11), (1.10) can be written as
$$H(\na\Phi)\bigl((\p_r\Phi)^2
+(\p_z\Phi)^2-q_0\p_z\Phi\bigr)-\rho_0 q_0\p_z\Phi+\rho_0
q_0^2=-\f{1}{\chi^2}H(\na\Phi)(\p_{\th}\Phi)^2\quad \text{on
$r=\chi(z,\th)$}.\tag3.4$$

Introducing the notation
$$\xi(z,\th)=\f{\chi(z,\th)-s_0 z}{z}.$$

Then (3.4) can be rewritten as
$$B_1\p_r\dot\vp+B_2\p_z\dot\vp
+B_3\xi=\kappa(\xi,\na_x\dot\vp)\quad \text{on $r=\chi(z,\th)$},
\tag3.5$$
where
$$
\left\{ \enspace
\aligned
B_1&=-\ds\f{\rho(s_0)}{c^2(\rho(s_0))}\biggl(u_{r}^2(s_0)
+u_z(s_0)(u_z(s_0)-q_0)\biggr)u_{r}(s_0)
+2\rho(s_0)u_{r}(s_0),\\
B_2&=-\ds\f{\rho(s_0)}{c^2(\rho(s_0))}\biggl(u_{r}^2(s_0)
+u_z(s_0)(u_z(s_0)-q_0)\biggr)u_z(s_0)
+2\rho(s_0)(u_z(s_0)-q_0)+(\rho(s_0)-\rho_0)q_0,\\
B_3&=\rho(s_0)\biggl(2u_{r}(s_0) u_r '(s_0)+2(u_z(s_0)-q_0)
{u_z}'(s_0)+q_0
{u_z}'(s_0)\biggr)+\rho'(s_0)\biggl(u_{r}^2(s_0)+{u_z}(s_0)({u_z}(s_0)-q_0)
\biggr)\\
&\qquad-\rho_0q_0 {u_z}'(s_0),
\endaligned
\right.
$$
and the generic function $\kappa(\xi,\nabla_x\dot\vp)$ is used to
denote a quantity dominated by
$C(b_0,q_0)|(\xi,\nabla_x\dot\vp)|^2$.

By Lemma 3.3 below, one has $B_1\not=0$ in (3.5) for large $q_0$.
Thus, Eq.~(3.5) can be rewritten as
$$\Cal B_0\dot\vp+\mu_2\xi=\kappa(\xi,\nabla_x\dot{\varphi})\quad
\text{on $r=\chi(z,\th)$};\tag3.6$$
here $\Cal B_0\dot\vp=\p_r\dot{\vp}+\mu_1\p_z\dot{\vp}$ with
$\mu_1=\ds\f{B_2}{B_1}$ and $\mu_2=\ds\f{B_3}{B_1}$.

Furthermore, (1.11) and the properties of $\Phi^-$ for $z\ge 1$ and
$\th\in [0,2\pi]$ listed in Lemma 3.1 imply
$$
\align \dot\vp(z, & \chi(z,\th),\th)=\Phi
(z,\chi(z,\th),\th)-{\hat\Phi}
(z,s_0z,\th)-\bigl({\hat\Phi}(z,\chi(z,\th),\th)-{\hat\Phi}
(z,s_0z,\th)\bigr)\\
&\quad =-\bigl({\hat\Phi}(z,\chi(z,\th),\th)-{\hat\Phi}
(z,s_0z,\th)\bigr)\\
&\quad =-\biggl(\int_0^{1}{\hat u}_r
(s_0+t\xi(z,\th))dt\biggr)z\xi(z,\th). \tag3.7
\endalign
$$

On the other hand, it follows from Lemma 3.1 that $\Phi^--q_0z\in
C_0^{\infty}(\O_-)$ in (1.7) and that, near the vertex of the cone
$r=b_0 z$, the solution $\Phi^+(x)$ is actually the background
solution. Thus, in order to prove Theorem 1.1, by the local existence
and stability result for the multi-dimensional shock solution to the
potential flow equation in [24], we only need to solve problem (3.1)
in the domain $\{(z,r,\th)\colon z\geq 1,\, b_0 z\leq
r\leq\chi(z,\th),\, 0\le\th\le 2\pi\}$ with the boundary conditions
(3.3), (3.6)--(3.7), and small initial data $\dot\vp(z,r,\th)|_{z=1}$,
$\p_z\dot\vp(z,r,\th)|_{z=1}$, and $\xi(z,\th)|_{z=1}$. Here,
smallness of initial data means that
$$\sum_{l\le k_0}|\na_x^l\dot\vp(1,r,\th)|
+\sum_{l\le k_0}|\na^l_{z,\th}\xi(1,\th)|\le C\ve,
\tag 3.8$$
where $k_0\in\Bbb N$, $k_0\ge 7$.

For later use, we now list specific estimates of the coefficients in
(3.1), (3.2), and (3.6). Since these estimates result from a direct,
but tedious computation that makes use of Lemmas 2.1 and 2.2, we
postpone the proof to the appendix.

\medskip

With respect to the coefficients of $\Cal L\dot\vp$ in (3.2), one has:

\smallskip

{\bf Lemma 3.2.} {\it For $q_0$ large enough, $1<\gamma<3$,
$0<b_0<b_*$, and $b_0\le s\le s_0$, one has
$$
\align
&P_1(s)=\f{b_0}{1-\f{\gamma-1}{2}b_0^2(1+b_0^2)}\biggl(1+O(\c)
+O(\b)\biggr),\\
&P_2(s)=\f{b_0^2\bigl(\f{3-\g}{2}-\f{\gamma-1}{2}b_0^2\bigr)}{1
-\f{\gamma-1}{2}b_0^2(1+b_0^2)}\biggl(1+O(\c)+O(\b)\biggr),\\
&P_3(s)=\f{\f{\g-1}{2}b_0^2(1+b_0^2)}{1-\f{\g-1}{2}b_0^2(1+b_0^2)}
\biggl(1+O(\c)+O(\b)\biggr),\\
&P_4(s)=O(\c)+O(\b),\\
&P_5(s)=-\f{\f{\g-1}{4}b_0^2(1+b_0^2)}{1-\f{\g-1}{2}b_0^2(1+b_0^2)}
\biggl(1+O(\c)+O(\b)\biggr),\\
&P_1'(s)=\f{-1+\f{\g-3}{2}b_0^2-\f{\g-1}{2}b_0^6}{(1+b_0^2)
\bigl(1-\f{\g-1}{2}b_0^2(1+b_0^2)\bigr)^2}
\biggl(1+O(\c)+O(\b)\biggr),\\
&P_2'(s)=\f{-2b_0+2(\g-2)b_0^3+2(\g-1)b_0^5}{(1+b_0^2)
\bigl(1-\f{\g-1}{2}b_0^2(1+b_0^2)\bigr)^2}
\biggl(1+O(\c)+O(\b)\biggr),\\
&P_3'(s)=-\f{(\g-1)b_0^3}{\bigl(1-\f{\g-1}{2}b_0^2(1+b_0^2)\bigr)^2}\biggl(1+
O(\c)+\f{1}{b_0^2} O(\b)\biggr).
\endalign$$}

\smallskip

Moreover, $B_i$, $i=1,2,3)$, in (3.5) and $\mu_j$, $j=1,2$, in (3.6)
admit the following estimates:

\smallskip

{\bf Lemma 3.3.}  {\it For large $q_0$, one has
$$
\align &B_1=\f{2}{1+b_0^2}\biggl(\f{(\gamma-1)}{2A\gamma(1
+b_0^2)}\biggr)^{\f{1}{\gamma-1}}(b_0 q_0)^{\f{\gamma+1}{\gamma-1}}
\biggl(1+O(\c)+O(\b)\biggr)>0,\\
&B_2=\f{1-b_0^2}{b_0 (1+b_0^2)}\biggl(\f{(\gamma-1)}{2A\gamma(1
+b_0^2)}\biggr)^{\f{1}{\gamma-1}}(b_0 q_0)^{\f{\gamma+1}{\gamma-1}}
\biggl(1+O(\c)+O(\b)\biggr),\\
&B_3=-\f{1}{b_0(1+b_0^2)^2}\biggl(\f{(\gamma-1)}{2A\gamma(1
+b_0^2)}\biggr)^{\f{1}{\gamma-1}}(b_0 q_0)^{\f{2\gamma}{\gamma-1}}
\biggl(1+\f{1}{b_0}O(\c)
+\f{1}{b_0}O(\b)\biggr),\\
&\mu_1=\f{1-b_0^2}{2b_0}\biggl(1+O(\c)+O(\b)\biggr)>0,\\
&\mu_2=-\f{q_0}{2(1+b_0^2)}\biggl(1+\f{1}{b_0}O(\c)+\f{1}{b_0}O(\b)\biggr)<0.
\endalign
$$}

\smallskip

{\bf Remark 3.1.} {\it We emphasize that there is a large factor 
$q_0$ in the coefficient $\mu_2$ of the boundary
condition\/ \rom{(3.6)} which shows that one has to be very careful in
treating the shock boundary condition later on.}

\vskip 0.4 true cm \centerline{\bf \S4. A first-order
weighted energy estimate} \vskip 0.3 true cm

In this section, we establish a weighted energy estimate of
$\na_x\dot\vp$ for the linear part of (3.1), together with (3.3) and
(3.6)--(3.8), which will play a fundamental role in our subsequent
analysis.

Set $D_{T}=\{(z,r,\th)\colon 1\le z\le T, \, b_0 z\le
r \le\chi(z,\th), \, 0\le\th <2\pi\}$ for any $T>1$.
$\Gamma_{T}=\{(z,r,\th)\colon r=\chi(z,\th), \, 1\le z\le T, \,
0\le\th <2\pi\}$ and $B_{T}=\{(z,r,\th)\colon r=b_0 z, 1\le z\le T,
0\le\th <2\pi\}$ are the lateral boundaries of $D_T$.

\medskip

{\bf Theorem 4.1.} {\it Let $\dot\vp\in
C^{\infty}(D_{T})$ satisfy the boundary condition\/ \rom{(3.3)}. Let
$|\xi(z,\th)|+|z\p_z\xi(z,\th)|+|\p_{\th}\xi(z,\th)|+|\na_x\dot\vp|\le
M\ve$ hold true for small $\ve$, $z\in [1, T]$, and some positive
constant $M$. Then there exists a multiplier $\Cal
M\dot\vp=A(z,r)\p_{z}\dot\vp+B(z,r)\p_{r}\dot\vp$ with smooth
coefficients such that, for any fixed constant $\mu<-1$,
$$ \align C_1  T^{\mu+1} & \iint_{b_0T\le
r\le\chi(T,\th)}(\na_x\dot\vp)^2(T,r,\th)
drd\th+C_2\iiint_{D_T}{z}^{\mu}(\na_x\dot\vp)^2 \,dzdrd\th\\
&\qquad +C_3\iint_{\Gamma_T}z^{\mu+1}(\p_z\dot\vp)^2
\,dS+C_4\iint_{\Gamma_{T}}{z}^{\mu+1}\f{(\p_{\th}\dot\vp)^2}{r^2}
\,dS\\
&\le \iiint_{D_T}\Cal L\dot\vp\cdot\Cal M\dot\vp dzdrd\th+C_5
\iint_{b_0 \leq
r\leq\chi(1,\th)}\biggl((\p_{z}\dot\vp)^2+(\p_{r}\dot\vp)^2
+\f{(\p_{\th}\dot\vp)^2}{r^2}\biggr)
(1,r,\th)\,drd\th\\
&\qquad +C_6\iint_{\Gamma_T}{z}^{\mu+1}(\Cal B_0\dot\vp)^2\,dS,\tag
4.1\endalign
$$
where $C_i$, $1\le i\le 6$, are some positive constants depending on
$b_0$ and $q_0$. In particular,
$$\cases
C_3&=\ds\f{(\g-1)b_0^2(1+b_0^2)^3}{8\bigl(1-\f{\g-1}{2}b_0^2(1+b_0^2)\bigr)}
+\f{1}{b_0^2}O((b_0
q_0)^{-2}) +\f{1}{b_0^2}O((b_0 q_0)^{-\f{2}{\g-1}}),\\
C_6&=\ds\f{(\g-1)b_0^4(1+b_0^2)}{2\bigl(1-\f{\g-1}{2}b_0^2(1+b_0^2)\bigr)}
+O((b_0 q_0)^{-2})+O((b_0 q_0)^{-\f{2}{\g-1}}).
\endcases\tag 4.2$$}

\smallskip

{\bf Remark 4.1.} {\it The values of constants $C_3$ and $C_6$ will
play an important role in the energy estimates for the nonlinear
problem\/ \rom{(3.3)} with \rom{(3.6)--(3.8)}. The most troublesome
term $\iint_{\Gamma_T}{z}^{\mu+1}(\Cal B_0\dot\vp)^2\,dS$ in the
right-hand side of\/ \rom{(4.1)} will be shown to be absorbed by
positive integrals in the left-hand side of\/ \rom{(4.1)}. The reason
for which the term $\iint_{\Gamma_T}{z}^{\mu+1}(\Cal
B_0\dot\vp)^2\,dS$ is hard is the following one: In view of the
Neumann boundary condition\/ \rom{(3.3)}, which is different from the
artificial Dirichlet boundary condition used in\/ \rom{[31]}, the
usual Poincar\'e inequality is not available for the solution
$\dot\vp$ on the shock surface which means that the
$L^2(\Gamma_T)$-estimates of $\na_x\dot\vp$ cannot be obtained
directly. \rom(Note that the boundary condition\/ \rom{(3.6)} contains
the function $\xi$, and that $\xi$ is roughly equivalent to
$\ds\f{\dot\vp}{z}$ in view of\/ \rom{(3.7)}, so that an estimate of
$\dot\vp$ on the shock surface has to be established.\rom)}

\smallskip

{\bf Proof.} We will determine the coefficients
$A(z,r)=z^{\mu}ra(\f{r}{z})$ and $B(z,r)=z^{\mu+1}b(\f{r}{z})$ for
$z\ge 1$ later on.  Set $\Cal
M\dot\vp=A(z,r)\p_{z}\dot\vp+B(z,r)\p_{r}\dot\vp$. Then it follows
from an integration by parts that
$$
\align &\iiint_{D_T} \Cal L\dot\vp\cdot\Cal M\dot\vp\, dzdrd\th
=\iiint_{D_T}z^{\mu}\biggl(K_1
(\p_{z}\dot\vp)^2+K_2\p_{z}\dot\vp\p_{r}\dot\vp
+K_3(\p_r\dot\vp)^2+K_4\f{(\p_{\th}\dot\vp)^2}{r^2}\biggr)\,dzdrd\th\\
&\qquad\quad +T^{\mu+1}\iint_{b_0 T\leq
r\leq\chi(T,\th)}N_1(T,r,\th)\, dS
-\iint_{b_0\leq r\leq \chi(1,\th)}N_1(1,r,\th)\, dS\\
&\qquad\quad +\iint_{\Gamma_T}z^{\mu+1}\biggl(N_2-\p_{z}\chi
N_1-\p_{\th}\chi N_3\biggr)\,dS+\iint_{B_T}z^{\mu+1}\biggl(b_0
N_1-N_2\biggr)\,dS,\tag 4.3
\endalign
$$
where
$$
\cases
&K_1=\biggl(\f{1}{2}s^2-sP_1(s)\biggr)a'(s)+\f{1}{2}b'(s)+\biggl(-\f{\mu}{2}s
-P_1(s)-sP_1'(s)+2P_4(s)\biggr)a(s),\\
&K_2=-sP_2(s)a'(s)+sb'(s)+\biggl(-P_2(s)+2P_5(s)-sP_2'(s)\biggr)a(s)
+\biggl(-(\mu+1)
+\f{2}{s}P_4(s)\biggr)b(s),\\
&K_3=-\f{1}{2}s^2
P_2(s)a'(s)+\biggl(sP_1(s)-\f{1}{2}P_2(s)\biggr)b'(s)
+\biggl(\f{\mu}{2}sP_2(s)
-\f{1}{2}s^2 P_2'(s)\biggr)a(s)\\
&\qquad\quad+\biggl(-(\mu+1)P_1(s)+sP_1'(s)-\f{1}{2}P_2'(s)+\f{2}{s}P_5(s)
\biggr)b(s),\\
&K_4=\f{1}{2}s^2 P_3(s) a'(s)-\f{1}{2}P_3(s)
b'(s)+\biggl(-\f{\mu}{2}sP_3(s)
+\f{1}{2}s^2 P_3'(s)
\biggr)a(s)+\biggl(\f{1}{s}P_3(s)-\f{1}{2}P_3'(s)\biggr)b(s)
\endcases\tag4.4
$$
and
$$
\cases
&N_1=\f{1}{2}sa(s)(\p_{z}\dot\vp)^2+b(s)\p_{z}\dot\vp\p_{r}\dot\vp
+\biggl(b(s)P_1(s)-\f{1}{2}sa(s)
P_2(s)\biggr)(\p_{r}\dot\vp)^2+sa(s) P_3(s)\ds\f{(\p_{\th}\dot\vp)^2}{2r^2},\\
&N_2=\biggl(sa(s)P_1(s)-\f{1}{2}b(s)\biggr)(\p_{z}\dot\vp)^2+sa(s)P_2(s)
\p_{z}\dot\vp\p_{r}\dot\vp
+\f{1}{2}b(s)P_2(s)(\p_r\dot\vp)^2+b(s)P_3(s)\ds
\f{(\p_{\th}\dot\vp)^2}{2r^2},\\
&N_3=-\f{1}{r^2}sa(s)P_3(s)\p_{z}\dot\vp\p_{\th}\dot\vp-\f{1}{r^2}b(s)P_3(s)
\p_{r}\dot\vp\p_{\th}\dot\vp.
\endcases\tag 4.5
$$

Our goal is to choose positive functions $a(s)$ and $b(s)$ in such a
way that all integrals over $D_T$, $B_T$, and $z=T$ in the right-hand
side of (4.3) are either positive or zero and that the integral over
$\Gamma_T$ provides sufficient control on $\dot\vp$. Starting from
here, we will derive sufficient conditions on the choice of $a(s)$ and
$b(s)$ in the process of analyzing each integral from which $a(s)$ and
$b(s)$ can then be determined. This process is subdivided into the
following five steps.

\vskip 0.2 true cm

{\bf Step 1.}  {\bf Handling of $\displaystyle\iint_{B_T}
z^{\mu+1}\biggl(b_0 N_1-N_2\biggr)dS$.}
By the boundary condition (3.3) and $u_r(b_0)=b_0 u_z(b_0)$, one
obtains that
$$\align
b_0
N_1-N_2&=\biggl(\f{1}{2}(1+b_0^2)+\f{b_0^2(1+b_0^2)u_z^2(b_0)}{2(u_z^2(b_0)
-c^2(\rho(b_0)))}\biggr)(b(b_0)-b_0^2
a(b_0))(\p_{z}\dot\vp)^2\\
&\qquad -(b(b_0) -b_0^2
a(b_0))P_3(b_0)\f{(\p_{\th}\dot\vp)^2}{2r^2}\endalign$$
on $B_T$. Since the terms
$\displaystyle\f{b_0^2(1+b_0^2)u_z^2(b_0)}{2(u_z^2(b_0)-c^2(\rho(b_0)))}$
and $P_3(b_0)$ are positive according to (iv) and (vii) in Lemma 2.1 and Lemma 3.2, in order to
control the integral $\displaystyle\iint_{B_T} z^{\mu+1}\biggl(b_0
N_1-N_2\biggr)dS$, one should choose
$$b(b_0)=b_0^2 a(b_0).\tag 4.6$$
This then implies
$$\iint_{B_T} z^{\mu+1}\biggl(b_0
N_1-N_2\biggr)\,dS=0.\tag 4.7$$

\vskip 0.2 true cm

{\bf Step 2.} {\bf Positivity of $\ds\iint_{b_0 T\leq
r\leq\chi(T,\th)}N_1 dS$.}
To ensure the positivity of the terms in $N_1$, which are quadratic in
$(\p_{z}\dot\vp, \p_r\dot\vp, \ds\f{\p_{\th}\dot\vp}{r})$, $a(s)$ and
$b(s)$ should fulfill
$$
\cases &a(s)>0,\\
&b^2(s)-2sP_1(s)a(s)b(s)+s^2P_2(s)a^2(s)<0
\endcases$$
which is equivalent to
$$a(s)>0,\quad \lambda_1(s)<\f{b(s)}{sa(s)}<\lambda_2(s).\tag 4.8$$
In this case, one arrives at
$$\iint_{b_0T\le r\le \chi(T,\th)} N_1\, dS\ge \iint_{b_0T\le
r\le\chi(T,\th)}\biggl(\lambda_{\min}(s)\bigl((\p_{z}\dot\vp)^2+(\p_r
\dot\vp)^2\bigr)+\f{1}{2r^2}saP_3(s)(\p_{\th}\dot\vp)^2\biggr)
\,dS,\tag 4.9
$$
where
$$\lambda_{\min}(s)=\f{1}{2}\biggl(\f{1}{2}sa(s)+b(s)P_1(s)-\f{1}{2}sa(s)P_2(s)
-\sqrt{\bigl(\f{1}{2}sa(s)-b(s)P_1(s)+\f{1}{2}sa(s)P_2(s)\bigr)^2+b(s)^2}
\biggr).$$
Moreover, by Lemmas 2.1, 2.2, and 3.2, the assumption on $\xi(z,\th)$
in Theorem 4.1, and the choices of $a(s)$, $b(s)$ to be made later on,
one actually obtains that
$$\lambda_{\min}(s)\ge C(b_0,q_0)+O(\ve).\tag 4.10$$

\vskip 0.3 true cm

{\bf Step 3.} {\bf Positivity of the integral over $D_T$.}
Under the constraints (4.6) and (4.8), we will choose $a(s)$ and
$b(s)$ in such a way that
$$K_1 (\p_{z}\dot\vp)^2+K_2\p_{z}\dot\vp\p_{r}\dot\vp+K_3(\p_r\dot\vp)^2
+K_4\f{(\p_{\th}\dot\vp)^2}{r^2}\ge 0.$$
This will be true if the coefficients $K_i$, $1\leq i\leq 4$, satisfy
$$\cases
&K_1>0,\\
&K_2^2-4K_1 K_3<0,\\
&K_4>0.
\endcases\tag4.11$$
According to (4.6), one can choose
$$\cases
&a(s)=1,\\
&b(s)=s^2\biggl(1+\ds\f{\t b}{b_0}(s-b_0)\biggr),
\endcases\tag 4.12$$
where the constant $\t b$ will  be determined in a moment.

It follows from Lemmas 2.1, 2.2, and 3.2 and a direct computation that
$$\cases
&K_1=\biggl(\ds\f{1}{2}b_0+O(b_0^2)\biggr)\t b-(\f{\mu}{2}-1)b_0
+O(b_0^2)+O((b_0q_0)^{-2})
+O((b_0q_0)^{-\f{2}{\g-1}}),\\
&K_2=\biggl(b_0^2+O(b_0^3)\biggr)\t
b+\bigl(2-\mu\bigr)b_0^2+O(b_0^3)+O((b_0q_0)^{-2})
+O((b_0q_0)^{-\f{2}{\g-1}}),\\
&K_3=\biggl(\ds\f{\g+1}{4}b_0^3+O(b_0^4)\biggr)\t
b+\biggl(1-\ds\f{\mu(\g+1)}{4}\biggr)b_0^3+O(b_0^4)+O((b_0q_0)^{-2})
+O((b_0q_0)^{-\f{2}{\g-1}}),\\
&K_4=-\biggl(\ds\f{\g-1}{4}b_0^3+O(b_0^4)\biggr)\t
b-\f{\mu}{4}(\g-1)b_0^3+O(b_0^4)+O((b_0q_0)^{-2})
+O((b_0q_0)^{-\f{2}{\g-1}}),\\
&K_2^2-4K_1K_3=-\ds\f{\g-1}{2}(\mu-\t b)(\mu-2-\t
b)b_0^4+O(b_0^5)+O((b_0q_0)^{-2}) +O((b_0q_0)^{-\f{2}{\g-1}}).
\endcases\tag 4.13$$
Therefore, one obtains from (4.11) that, for $q_0$ large enough,
$$\t b+2-\mu>0,\quad (\mu-\t b)(\mu-2-\t b)>0,\quad \t b+\mu<0.$$
If one sets
$$\t b=1,\tag 4.14$$
then
$$\mu<-1.\tag 4.15$$

In this case, one arrives at
$$\align \iiint_{D_T} {z}^{\mu} & \biggl\{K_1
(\p_{z}\dot\vp)^2+K_2\p_{z}\dot\vp\p_{r}\dot\vp+K_3(\p_r\dot
\vp)^2+K_4 (\p_{\th}\dot\vp)^2\biggr\}\,dzdrd\th\\
&\qquad \geq  C(b_0)\iiint_{D_T}z^{\mu}\biggl((\p_{z}\dot\vp)^2+
(\p_{r}\dot\vp)^2+\f{(\p_{\th}\dot\vp)^2}{r^2}\biggr)\,dzdrd\th.\tag
4.16\endalign$$

\vskip 0.2 true cm{\bf Step 4. Estimate of
$\displaystyle\iint_{\Gamma_T}z^{\mu+1} \biggl(N_2-\p_{z}\chi
N_1-\p_{\th}\chi N_3\biggr)\,dS$.}
By the assumptions on $\xi(z,\th)$ in Theorem 4.1, it follows from
the expressions for $N_1$, $N_2$, $N_3$ and a direct computation that
$$
N_2-\p_{z}\chi N_1-\p_{\th}\chi
N_3=\beta_0(\p_z\dot\vp)^2+\beta_1\p_z\dot\vp\p_r\dot\vp
+\beta_2(\p_r\dot\vp)^2+\beta_3\f{1}{r^2}(\p_{\th}\dot\vp)^2,\tag4.17
$$
where
$$\align
&\beta_0=s_0P_1(s_0)a(s_0)-\f{1}{2}s_0^2a(s_0)-\f12b(s_0)+O(\ve),\\
&\beta_1=s_0P_2(s_0)a(s_0)-s_0b(s_0)+O(\ve),\\
&\beta_2=\f12P_2(s_0)b(s_0)-s_0P_1(s_0)b(s_0)+\f12s_0^2P_2(s_0)a(s_0)+O(\ve),\\
&\beta_3=\f12P_3(s_0)(b(s_0)-s_0^2a(s_0))+O(\ve).
\endalign
$$
Due to $\p_r\dot\vp={\Cal B}_0\dot\vp-\mu_1\p_z\dot\vp$, from (4.17)
one obtains that
$$
N_2-\p_{z}\chi N_1-\p_{\th}\chi
N_3=\bigl(\beta_0-\mu_1\beta_1+\beta_2\mu_1^2\bigr)(\p_z\dot\vp)^2+
\bigl(\beta_1-2\mu_1\beta_2\bigr)\p_z\dot\vp{\Cal B}_0\dot\vp
+\beta_2({\Cal
B}_0\dot\vp)^2+\beta_3\f{1}{r^2}(\p_{\th}\dot\vp)^2.
$$

By Lemmas 2.1 and 3.2 and the expressions for $a(s), b(s)$ in
(4.12) and (4.14), one arrives at
$$
\align
&\beta_0-\mu_1\beta_1+\beta_2\mu_1^2=\f{1}{4(1-\f{\g-1}{2}b_0^2(1+b_0^2))}
\biggl(\f{\g-1}{2}b_0^2(1+b_0^2)^3+O(\c)+O(\b)+O(\ve)\biggr),\\
&\beta_1-2\mu_1\beta_2=O(\c)+O(\b)+O(\ve),\\
&\beta_2=\f{1}{1-\f{\g-1}{2}b_0^2(1+b_0^2)}\biggl(-\f{\g-1}{2}b_0^4(1+b_0^2)
+O(\c)+O(\b)+O(\ve)\biggr),\\
&\beta_3\ge C(b_0,q_0)+O(\ve).
\endalign
$$

Consequently,
$$\align
\iint_{\Gamma_T}z^{\mu+1} &\biggl(N_2-\p_{z}\chi N_1-\p_{\th}\chi N_3\biggr)dS\\
&\geq
(Q_1(b_0,q_0)+O(\ve))\iint_{\Gamma_T}z^{\mu+1}(\p_z\dot\vp)^2
dS-(Q_2(b_0,q_0)+O(\ve))\iint_{\Gamma_T}z^{\mu+1}(\Cal
B_0\dot\vp)^2\,dS\\
&\qquad +\bigl(C(b_0,q_0)+O(\ve)\bigr)\iint_{\Gamma_T}z^{\mu+1}
\f{(\p_{\th}\dot\vp)^2}{r^2}\,dS,\tag4.18\endalign$$
where
$$
\align
&Q_1(b_0,q_0)=\f{(\g-1)b_0^2(1+b_0^2)^3}{8\bigl(1-\f{\g-1}{2}
b_0^2(1+b_0^2)\bigr)}+O(\c)+O(\b)>0,\\
&Q_2(b_0,q_0)=\f{(\g-1)b_0^4(1+b_0^2)}{2\bigl(
1-\f{\g-1}{2}b_0^2(1+b_0^2)\bigr)}+O(\c)+O(\b)>0.
\endalign
$$

\vskip 0.2 true cm {\bf Step 5. Estimate of
$\displaystyle\iint_{b_0 \leq r\leq
\chi(1,\th)}N_1(1,r,\th)\,dS$.}
>From the expression for $N_1$ and the initial condition (3.8), one
easily obtains
$$\biggl|\iint_{b_0\leq r\leq
\chi(1,\th)}K_5(1,r,\th)drd\th\biggr| \leq C(b_0,q_0)\ve^2.\tag
4.19$$

Finally, substituting the estimates (4.7), (4.9)--(4.10), (4.16), and
(4.18)--(4.19) into (4.3), (4.1) and (4.2) are obtained in terms of
the smallness of $\ve$ given in Remark 1.1. Therefore, Theorem 4.1 is
proved. \qed

\smallskip

Based on Theorem 4.1, we will derive a first-order uniform energy
estimate of $\na_x\dot\vp$ for the linear part of (3.1) together with
(3.3) and (3.6)--(3.8). To this end, we require a Hardy-type
inequality on $\iint_{\Gamma_T}z^{\mu-1}|\dot\vp|^2 dS$, which is
motivated by [14, Theorem~330] and derived utilizing the special
structures of (3.6)--(3.7).

\medskip

{\bf Lemma 4.2.(Hardy-type inequality)} {\it Under the assumptions
of\/ \rom{Theorem 4.1}, for $\mu<-1$, one has
$$ \align
&\iint_{\Gamma_T}z^{\mu-1}|\dot\vp|^2\, dS\le
C(b_0,q_0)\ve^2+\f{(1+b_0^2)^2}{\mu^2}\bigl(1+O(\c)+O(\b)\bigr)
\iint_{\Gamma_T}z^{\mu+1}(\p_z\dot\vp)^2\,dS\\
&\quad +C(b_0,q_0)\ve\biggl(\iint_{\Gamma_T}z^{\mu+1}(\Cal
B_0\dot\vp)^2\,dS
+\iint_{\Gamma_T}z^{\mu+1}\f{(\p_{\th}\dot\vp)^2}{r^2}\,dS\biggr).\tag4.20
\endalign
$$}

{\bf Proof.} First we assert that
$$
\align \iint_{\Gamma_T}z^{\mu-1}|\dot\vp|^2 dS &\le
\bigl(\f{4}{\mu^2}+O(\ve)\bigr)\iint_{\Gamma_T}
z^{\mu+1}(1-\mu_1\p_z\chi)^2(\p_z\dot\vp)^2\,dS \\
& \qquad +C(b_0,q_0)\ve\biggl(\iint_{\Gamma_T}z^{\mu+1}({\Cal
B_0}\dot\vp)^2\,dS
+\iint_{\Gamma_T}z^{\mu+1}\f{(\p_{\th}\dot\vp)^2}{r^2}\,dS\biggr)
+C(b_0,q_0)\ve^2.\tag
4.21
\endalign$$

Indeed, note that
$$\iint_{\Gamma_T}z^{\mu-1}\dot\vp^2 dS
=\int_0^{2\pi}d\th\int_1^Tz^{\mu-1}\dot\vp^2(z,\chi(z,\th),\th)dz.\tag4.22$$
Set $m(\th)\equiv
\int_1^Tz^{\mu-1}\dot\vp^2(z,\chi(z,\th),\th)dz$. Then, by an
integration by parts,
$$
\align
m(\th)&=\f{1}{\mu}z^{\mu}\dot\vp^2(z,\chi(z,\th),\th)\biggr|_{z=1}^{z=T}
-\f{2}{\mu}\int_1^Tz^{\mu}\dot\vp(z,\chi(z,\th),\th)(\p_z\dot\vp
+\p_z\chi\p_r\dot\vp)(z,\chi(z,\th),\th)dz\\
&\le\f{1}{|\mu|}\dot\vp^2(1,\chi(1,\th),\th)-\f{2}{\mu}\int_1^T
z^{\mu}\dot\vp(z,\chi(z,\th),\th)\bigl(\p_z\chi\Cal
B_0\dot\vp+(1-\mu_1\p_z\chi)\p_z\dot\vp\bigr)(z,\chi(z,\th),\th)dz.\tag4.23
\endalign$$

By (3.7), one has
$$\dot\vp(z,\chi(z,\th),\th)=a_1(z,\chi(z,\th),\th)\xi,\tag 4.24$$
where
$$a_1(z,\chi(z,\th),\th)\equiv -\biggl(\int_0^{1}{\hat u}_r
(s_0+t\xi(z,\th))dt\biggr)z<0.\tag 4.25$$

Due to the assumptions in Theorem 4.1 and Lemma 2.1,
$\p_z\chi=b_0(1+O(\b)+O(\ve))$, and thus because of the smallness of $\ve$,
$$1-\mu_1\p_z\chi=\f{1+b_0^2}{2}(1+O(\b)>0.\tag4.26$$

It follows from (3.8), (3.6), $\mu_2<0$ by Lemma 3.3, (4.24)--(4.25),
and a direct computation that
$$
\align m(\th)&\le
C(b_0,q_0)\ve^2-\f{2}{\mu}\int_1^Tz^{\mu}\dot\vp(z,\chi(z,\th),\th)
(1-\mu_1\p_z\chi)
\p_z\dot\vp(z,\chi(z,\th),\th)\,dz\\
&\qquad
+\f{2\mu_2}{\mu}\int_1^Tz^{\mu}\p_z\chi\dot\vp(z,\chi(z,\th),\th)\xi
\,dz-\f{2}{\mu}\int_1^Tz^{\mu}\p_z\chi\dot\vp(z,\chi(z,\th),\th)\kappa(\xi,
\na\dot\vp)\,dz\\
&=C(b_0,q_0)\ve^2-\f{2}{\mu}\int_1^Tz^{\mu}\dot\vp(z,\chi(z,\th),\th)
(1-\mu_1\p_z\chi)\p_z\dot\vp(z,\chi(z,\th),\th)\,dz\\
&\qquad +\f{2\mu_2}{\mu}\int_1^Tz^{\mu}\p_z\chi
a_1(z,\chi(z,\th),\th)\xi^2
\,dz-\f{2}{\mu}\int_1^Tz^{\mu}\p_z\chi\dot\vp(z,\chi(z,\th),\th)\kappa(\xi,
\na\dot\vp)\,dz\\
&\le
C(b_0,q_0)\ve^2+\f{1}{2}\int_1^Tz^{\mu-1}\dot\vp^2(z,\chi(z,\th),\th)\,dz
+\f{2}{\mu^2}\int_1^Tz^{\mu+1}(1-\mu_1\p_z\chi)^2
(\p_z\dot\vp)^2(z,\chi(z,\th),\th)\,dz\\
&\qquad
+C(b_0,q_0)\ve\biggl(\int_1^Tz^{\mu-1}\dot\vp^2(z,\chi(z,\th),\th)\,dz+
\int_1^Tz^{\mu+1}\xi^2\,dz+\int_1^Tz^{\mu+1}\bigl(|\na_{r,z}\dot\vp|^2\\
&\qquad
+\f{(\p_{\th}\dot\vp)^2}{r^2}\bigr)(z,\chi(z,\th),\th)\,dz\biggr).
\endalign
$$
Together with (4.22) and (4.26), this yields
$$
\align \iint_{\Gamma_T}& z^{\mu-1} \dot\vp^2 \,dS\\
&\le C(b_0,q_0)\ve^2
+(\f{4}{\mu^2}+O(\ve))\iint_{\Gamma_T}z^{\mu+1}(1-\mu_1\p_z\chi)(\p_z\dot\vp)^2
\,dS\\
&\qquad +C(b_0,q_0)\ve\biggl(\iint_{\Gamma_T}z^{\mu+1}(\Cal
B_0\dot\vp)^2
\,dS+\iint_{\Gamma_T}z^{\mu+1}\f{(\p_{\th}\dot\vp)^2}{r^2}\,dS\biggr)\\
&\le
C(b_0,q_0)\ve^2+\f{(1+b_0^2)^2}{\mu^2}\bigl(1+O(\c)+O(\b)\bigr)
\iint_{\Gamma_T}z^{\mu+1}(\p_z\dot\vp)^2\,dS\\
&\qquad +C(b_0,q_0)\ve\biggl(\iint_{\Gamma_T}z^{\mu+1}(\Cal
B_0\dot\vp)^2
\,dS+\iint_{\Gamma_T}z^{\mu+1}\f{(\p_{\th}\dot\vp)^2}{r^2}\,dS\biggr).
\endalign
$$
Hence, Lemma 4.2 is proved. \qed

\medskip

{\bf Theorem 4.3.} {\it Under the assumptions of\/ \rom{Theorem 4.1}, for
$\mu<-1$, one has
$$ \align C_0T^{\mu+1} & \iint_{b_0T\le
r\le\chi(T,\th)}|\na_x\dot\vp(T,r,\th)|^2 \,drd\th
+C_0\iiint_{D_T}{z}^{\mu}|\na_x\dot\vp|^2 \, dzdrd\th\\
&\qquad+C_0\iint_{\Gamma_{T}}{z}^{\mu+1}|\na_x\dot\vp|^2\, dS
\leq\iiint_{D_T}\Cal L\dot\vp\cdot\Cal M\dot\vp
\,dzdrd\th+C(b_0,q_0)\ve^2,\tag 4.27\endalign
$$
where $C_0=C_0(b_0, q_0)>0$ is a generic constant.}

{\bf Proof.} To obtain (4.27), we are required to give a
delicate estimate of the term
$C_6\displaystyle\iint_{\Gamma_T}z^{\mu+1}(\Cal B_0\dot\vp)^2\, dS$ in
the right-hand side of (4.1) so that it can be absorbed by the
positive terms in the left-hand side of (4.1).

We now treat the term $\displaystyle\iint_{\Gamma_T}z^{\mu+1}(\Cal
B_0\dot\vp)^2\, dS$.

>From (3.6) and the definition of $\kappa(\xi,\na_x\dot\vp)$, one has
$$\align \ds\iint_{\Gamma_T}z^{\mu+1} & (\Cal
B_0\dot\vp)^2\,dS =\ds\iint_{\Gamma_T}z^{\mu+1}(\kappa(\xi,\nabla_{x}\dot\vp)
-\mu_2\xi)^2\,dS\\
&\leq\mu_2^2(1+O((b_0 q_0)^{-2}))\ds\iint_{\Gamma_T}
z^{\mu-1}(z\xi)^2\,dS+C(b_0,q_0)\iint_{\Gamma_T}
z^{\mu+1}\kappa^2(\xi,\nabla_{x}\dot\vp)\,dS.\tag4.28
\endalign
$$
Note that
$$
\align \ds\iint_{\Gamma_T}
z^{\mu+1}\kappa^2(\xi,\nabla_{x}\dot\vp)\,dS &\ds\leq
C(b_0,q_0)\ve^2\iint_{\Gamma_T}z^{\mu+1}\bigl(\xi^2
+|\nabla_x\dot\vp|^2\bigr)\,dS\\ &\leq\ds
C(b_0,q_0)\ve^2\iint_{\Gamma_T}z^{\mu+1}\biggl(\xi^2+(\p_{z}\dot\vp)^2
+\f{(\p_{\th}\dot\vp)^2}{r^2}+(\Cal
B_0\dot\vp)^2\biggr)\,dS.\tag4.29\endalign
$$
Furthermore, as $\ve$ is small and by the boundary condition (3.7)
together with Lemma 2.1 (iii), one obtains from (4.28)--(4.29) that
$$\align
\iint_{\Gamma_T}z^{\mu+1}(\Cal B_0\dot\vp)^2 \,dS
&\leq  \f{\mu_2^2(1+b_0^2)^2}{(b_0
q_0)^2}\biggl(1+O(\c)+O(\b)\biggr)\iint_{\Gamma_T}z^{\mu-1}\dot\vp^2
\,dS\\
&\qquad +C(b_0,q_0)\ve^2\iint_{\Gamma_T}z^{\mu+1}\biggl((\p_z\dot\vp)^2
+\f{(\p_{\th}\dot\vp)^2}{r^2}
\biggr)\,dS\\
&\leq \f{1}{4b_0^2}\biggl(1+\f{1}{b_0}O(\c)+\f{1}{b_0}O(\b)
\biggr)\iint_{\Gamma_T}z^{\mu-1}\dot\vp^2
\,dS\\
&\qquad +C(b_0,q_0)\ve^2\iint_{\Gamma_T}z^{\mu+1}\biggl((\p_z\dot\vp)^2
+\f{(\p_{\th}\dot\vp)^2}{r^2}\biggr)\,dS.\tag 4.30
\endalign$$
Therefore, combining  (4.30) with (4.20) in Lemma 4.2 yields
$$\align
\iint_{\Gamma_T}z^{\mu+1} & (\Cal B_0\dot\vp)^2\, dS\\
& \leq
C(b_0,q_0)\ve^2+\f{(1+b_0^2)^2}{4b_0^2\mu^2}\biggl(1+\f{1}{b_0}O(\c)
+\f{1}{b_0}O(\b)
\biggr)\iint_{\Gamma_T}z^{\mu+1}(\p_z\dot\vp)^2 \,dS\\
&\qquad
+C(b_0,q_0)\ve\iint_{\Gamma_T}z^{\mu+1}\f{(\p_{\th}\dot\vp)^2}{r^2}\,dS.\tag
4.31
\endalign$$
Substituting (4.31) into (4.1), one obtains
$$ \align C_0  T^{\mu+1}&\iint_{b_0T\le
r\le\chi(T,\th)}|\na_x\dot\vp|^2 \,drd\th+C_0\iiint_{D_T}{z}^{\mu}
|\na_x\dot\vp|^2\, dzdrd\th \\
&\qquad +\biggl(Q_0(b_0,q_0)+\f{1}{b_0^2}O(\c)
+\f{1}{b_0^2}O(\b)\biggr)
\iint_{\Gamma_T}z^{\mu+1}(\p_z\dot\vp)^2
\,dS\\
&\qquad +C_0\iint_{\Gamma_{T}}{z}^{\mu+1}\f{(\p_{\th}\dot\vp)^2}{r^2}\,dS\\
&\le \iiint_{D_T}\Cal L\dot\vp\cdot\Cal M\dot\vp
\,dx+C(b_0,q_0)\ve^2,\tag 4.32\endalign
$$
where
$Q_0(b_0,q_0)=\ds\f{(\g-1)b_0^2(1+b_0^2)^3}{8\bigl(1-\f{\g-1}{2}b_0^2
(1+b_0^2)\bigr)}\bigl(1-\f{1}{\mu^2}\bigr)$ and $C_0=C_0(b_0, q_0)>0$
is a generic constant. Moreover,
$$Q_0(b_0,q_0)>0
$$
because of $\mu^2>1$.  Together with (4.32), this finishes the proof of
Theorem 4.3. \qed

\vskip 0.4 true cm \centerline{\bf \S5. Higher-order weighted
energy estimates} \vskip 0.4 true cm

In this section, we derive a higher-order energy estimate so that one
can establish the decay properties of $\na_x\dot\vp$ and $\xi$ for large
$z$.

As in [20], we denote by
$$S=\{S_1, S_2\},\quad
\text{where}\quad S_1=z\p_z+r\p_r,\quad S_2=\p_{\th},\tag 5.1$$
certain Klainerman vector fields and by
$$S_{\Gamma}=\{S_{1\Gamma},S_{2\Gamma}\},\quad \text{where}\quad
S_{1\Gamma}=z\p_z+z\p_{z}\chi(z,\th)\p_r,\quad
S_{2\Gamma}=\p_{\th}+\p_{\th}\chi\p_{r},\tag 5.2
$$
vector fields which are tangential to $\Gamma$. We will use these
vector fields, which are tangential to the cone surface and are nearly
tangential to the shock surface, respectively, to act upon Eq.~(3.1) and
the boundary conditions (3.3) and (3.6). This allows us to raise the
order of the energy estimates by a commutator argument. Let us point
out that there is a difference to the usual commutator argument
inasmuch as the vector fields $S$ are only nearly tangential to the shock
front which causes certain error terms to occur in the estimates which
in turn is due to the perturbation of the shock surface with $r=s_0z$.
Furthermore, we cannot adapt the analysis of [31] as we have to deal
with a Neumann-type boundary condition on the fixed boundary, while
[31] treats an artificial Dirichlet-type boundary condition so that
there a Poincar\'e-type inequality (see also [11, Lemma 1]) as one of
the key elements of the analysis of [31] is available. However, by
making use of the delicate energy estimate established in \S4, we will
be able to drive the desired estimates.

\medskip

To prove Theorem 1.1, we require the following elementary estimate:

\smallskip

{\bf Lemma 5.1.}  {\it Let $\dot \vp$ be a $C^{k_0}(D_T)$
 solution to \rom{(3.1)}, where $k_0\in\Bbb N$, and
$$\ds\sum_{0\leq l_1+l_2\leq
[\f{k_0}{2}]+1}z^{l_1}|\p_{z}^{l_1}\p_{\th}^{l_2}\xi|+\dsize\sum_{0\le
l\le [\f{k_0}{2}]+1}z^l|\na_x^{l+1}\dot\vp|\leq
M\ve \quad \text{\rom{in} $D_T$,}
$$
where $M>0$ is some constant and $\ve$ is sufficiently small.
Then
$$\dsize C\dsize\sum_{0\le l\le k_0-1}|\na_x
S^l\dot\vp|\leq\sum_{0\le l\le k_0-1}x_3^l|\na_x^{l+1}\dot\vp| \le
C(b_0,q_0)\dsize\sum_{0\le l\le k_0-1}|\na_x S^l\dot\vp|\quad \text{\rom{in}
$D_T$.}\tag5.3$$}

{\bf Proof.}  The first inequality is immediate from the definition of
the Klainerman vector fields $S$ in (5.1).

We show the second inequality in (5.3). The case $k_0=1$ can be
verified directly. Cases when $k_0\geq 3$ can then be obtained by an
inductive argument. Thus, we only need to deal with the case $k_0=2$.

Because of
$$\align
&\p_r=\f{1}{r}\biggl(S_1-z\p_z\biggr),\\
&\p_{zr}^2=\f{1}{r}\biggl(\p_z S_1-\p_z-z\p_{z}^2\biggr),\\
&\p_{r}^2=-\f{1}{r}\p_r+\f{1}{r}\p_r S_1-\f{1}{r^2}z\p_z
S_1+\f{z^2}{r^2}\p_{z}^2,\tag5.4
\endalign$$
it follows from Eq.~(3.1) and the assumptions in Lemma 5.1 that
$$
\biggl|\biggl((1-f_1)-(2P_1(s)-f_2)\f{1}{s}+(P_2(s)-f_3)\f{1}{s^2}\biggr)
\p_z^2\dot\vp\biggr|
 \le C(b_0,q_0)\bigl(\f{1}{r}|\na_x
S_2\dot\vp|+\f{1}{r}|\na_x\dot\vp|\bigr).\tag 5.5
$$
Moreover, by the assumptions in the Lemma 5.1 and the expression
of $f_i(1\leq i\leq 3)$, one has
$$\ds\sum_{i=1}^{3}|f_i|\leq C(b_0,q_0)\ve.
$$
Combining this with $s^2-2P_1(s)s+P_2(s)<0$ from Lemma 2.2 and
(5.4)--(5.5) yields
$$r|(\p_{z}^2\dot\vp, \p_{zr}^2\dot\vp,\p_{r}^2\dot\vp)|\leq
C(b_0,q_0)(|\na_x S\dot\vp|+|\na_x\dot\vp|).
$$
Together with  $S_2=\p_{\th}$ and the change of coordinates (1.5),
this concludes the proof of Lemma 5.1. \qed

\smallskip

Now we begin to establish the higher-order energy estimates.

Since the vector fields $S$ are tangential to the boundary $r=b_0 z$,
$\p_r S^m\dot\vp=b_0\p_z S^m\dot\vp$ holds on $r=b_0z$ in view of the
boundary condition (3.3). Therefore, one can apply Theorem 4.1
directly to $S^m\dot\vp$ ($0\le m\le k_0-1$) to obtain:

\smallskip

{\bf Lemma 5.2.} {\it Let the assumptions of\/ \rom{Theorem 4.1} be
fulfilled. If $\dot\vp$ is a $C^{k_0}(\overline{D_T})$-solution to
problem \rom{(3.1)} with \rom{(3.3)} and \rom{(3.6)--(3.8)}, where
$k_0\geq 7$, then, for $0\leq m\leq k_0-1$ and $\mu<-1$,
$$ \align C_1  T^{\mu+1} & \iint_{b_0T\le
r\le\chi(T,\th)}|\na_x S^{m}\dot\vp|^2(T,r,\th)
\,drd\th+C_2\iiint_{D_T}{z}^{\mu}
|\na_x S^{m}\dot\vp|^2 \,dzdrd\th\\
&\qquad+C_3\iint_{\Gamma_T}z^{\mu+1}(\p_z S^{m}\dot\vp)^2
\,dS+C_4\iint_{\Gamma_{T}}{z}^{\mu+1}\f{1}{r^2}(\p_{\th}
S^{m}\dot\vp)^2\,dS\\
& \le \iiint_{D_T}\Cal L S^{m}\dot\vp\cdot\Cal M S^{m}\dot\vp
\,dzdrd\th+C_5 \iint_{b_0 \leq
r\leq\chi(1,\th)}\biggl((\p_{z}S^{m}\dot\vp)^2 \\
& \qquad +(\p_{r}S^{m}\dot\vp)^2
+\f{1}{r^2}(\p_{\th}S^{m}\dot\vp)^2\biggr)(1,r,\th)\,drd\th
+C_6\iint_{\Gamma_T}{z}^{\mu+1}(\Cal B_0 S^{m}\dot\vp)^2\,dS.\tag
5.6\endalign
$$
Here, the constants $C_i$, $1\le i\le 6$, are given in\/ \rom{Theorem
4.1}.} \qed

As in Theorem 4.3, one needs to control the term
$\displaystyle\iint_{\Gamma_T}z^{\mu+1}(\Cal B_0 S^m\dot\vp)^2\,dS$
in (5.6) to obtain related higher-order weighted energy
estimates. In addition, the term  $\iiint_{D_T}S^{m}\Cal L
\dot\vp\cdot\Cal M S^{m}\dot\vp\, dzdrd\th$ appearing in the right-hand side
of (5.6) has also to be payed attention to. This will produce some
boundary terms on the shock surface and conic surface, respectively,
by integration by parts.

\medskip

We now establish:

\smallskip

{\bf Theorem 5.3.} {\it Let $\dot\vp\in C^{k_0}(\overline{D_T})$ and
$\xi(z,\th)\in C^{k_0}([0, 2\pi]\times [1, T])$ be solutions
to \rom{(3.1)} with \rom{(3.3)} and \rom{(3.6)--(3.8)}, where $k_0\ge
7$, and
$$\ds\sum_{0\leq l_1+l_2\leq
[\f{k_0}{2}]+1}z^{l_1}|\p_{z}^{l_1}\p_{\th}^{l_2}\xi|+\dsize\sum_{0\le
l\le [\f{k_0}{2}]+1}z^l|\na_x^{l+1}\dot\vp|\leq M\ve,\tag5.7$$ where
$M>0$ is a constant. Then, for sufficiently small $\ve>0$ and $\mu<-1$,
$$
\align &\iint_{b_0T\le r\le\chi(T,\th)}\dsize\sum_{0\le l\le
k_0-1}T^{2l+\mu+1} |\na_x^{l+1}\dot\vp|^2(T,r,\th)\,dS
+\iiint_{D_T}\dsize\sum_{0\le l
\le k_0-1}z^{2l+\mu}|\na_x^{l+1}\dot\vp|^2\,dzdrd\th\\
&\qquad +\iint_{\Gamma_T}\dsize\sum_{0\le l\le k_0-1}z^{2l+\mu+1}
|\na_x^{l+1}\dot\vp|^2\,dS\le C(b_0,q_0)\ve^2.\tag5.8
\endalign
$$}

{\bf Proof.} Note that on the shock surface $\Gamma$
$$\Cal B_0 S^m\dot\vp=S^m_{\Gamma}\Cal B_0 \dot\vp+\Cal B_0
(S^m-S^m_{\Gamma})\dot\vp+[\Cal B_0, S^m_{\Gamma}]$$ and
$$S^m_{\Gamma}\dot\vp=-\bigl(\int_0^{1}{\hat u}_r(s_0+t\xi)dt\bigr)z
S^m_{\Gamma}\xi+\bigl[S^m_{\Gamma}, -\bigl(\int_0^{1}{\hat
u}_r(s_0+t\xi)dt\bigr)z\bigr]\xi;$$ here and below $\bigl[\cdot ,
\cdot\bigr]$ stands for the commutator.
Together with Lemma 5.2 and the initial condition (3.8), this yields
$$
\align &T^{\mu+1}\iint_{b_0T\le r\le\chi(T,\th)}|\na_x
S^m\dot\vp|^2(T,r,\th)\, drd\th +\iiint_{D_T}{z}^{\mu}|\na_x
S^m\dot\vp|^2\,dzdrd\th
+\iint_{\Gamma_{T}}{z}^{\mu+1}|\na_x S^m\dot\vp|^2\,dS\\
&\qquad \leq C(b_0,q_0)\ve^2+C(b_0,q_0)\dsize\sum_{i=1}^4I_k,\tag
5.9\endalign
$$
where
$$\align
&I_1=\iiint_{D_T}\Cal L S^m\dot\vp\cdot\Cal M S^m\dot\vp
\,dzdrd\th,\\
&I_2=\iint_{\Gamma_T}z^{\mu-1}\biggl|\bigl[S^m_{\Gamma},
\bigl(\int_{0}^{1}{\hat u}_r(s_0+t\xi)dt\bigr)z\bigr]\xi\biggr|^2\, dS,\\
&I_3=\iint_{\Gamma_T}z^{\mu+1}|\Cal
B_0(S^m-S^m_{\Gamma})\dot\vp|^2\, dS,\\
&I_4=\iint_{\Gamma_T}z^{\mu+1}|[\Cal B_0, S^m_{\Gamma}]\dot\vp|^2
\,dS.
\endalign
$$
To prove Theorem 5.3, we treat each $I_k$ in (5.9) separately. With
this aim, we divide the rest of the proof into four steps.

\vskip 0.2 true cm

{\bf Step 1. Estimate of $I_1$.}
%
First we derive an explicit representation of $\Cal L
S^m\dot\vp$.
Due to $S_1(\ds\f{r}{z})=S_2(\ds\f{r}{z})=0$ and
$S_1(\ds\f{1}{r})=-\ds\f{1}{r}, S_2(\ds\f{1}{r})=0$,
$$\Cal
L S_1\dot\vp=S_1\Cal L \dot\vp-2 \Cal L\dot\vp,\quad \Cal L
S_2\dot\vp=S_2\Cal L \dot\vp.$$
By induction, one obtains
$$\Cal
L S^m\dot\vp=S^m\Cal L \dot\vp+\dsize\sum_{0\le l\le
m-1}C_{lm}S^l\Cal L \dot\vp.\tag5.10$$
Hence, it follows from Eq.~(3.1) and (5.10) that
$$S^m\Cal L\dot\vp=I_{11}+I_{12}+I_{13},\tag 5.11$$
where
$$\align
I_{11}&=f_1(s,\na_x\dot\vp)\p_{z}^2 S^m\dot\vp+f_2(s,\na_x\dot\vp)
\p_{zr}^2 S^m\dot\vp
+f_3(s,\na_x\dot\vp)\p_{r}^2 S^m\dot\vp+\f{1}{r^2}f_4(s,\na_x\dot\vp)
\p_{\th}^2 S^m\dot\vp\\
&\quad+\f{1}{r}f_5(s,\na_x\dot\vp)\p_{z\th}^2 S^m\dot\vp+\f{1}{r}
f_6(s,\na_x\dot\vp)\p_{r\th}^2 S^m\dot\vp,\\
I_{12}&=f_1(s,\na_x\dot\vp)[S^m, \p_{z}^2]
\dot\vp+f_2(s,\na_x\dot\vp)[S^m, \p_{zr}^2]\dot\vp
+f_3(s,\na_x\dot\vp)[S^m, \p_{r}^2]\dot\vp+\f{1}{r^2}f_4(s,\na_x\dot\vp)
[S^m, \p_{\th}^2]\dot\vp\\
&\quad+\f{1}{r}f_5(s,\na_x\dot\vp)[S^m, \p_{z\th}^2
]\dot\vp+\f{1}{r}f_6(s,\na_x\dot\vp)[S^m, \p_{r\th}^2]\dot\vp,
\endalign
$$
$$
\align
I_{13}&=\dsize\sum_{0\le l\le
m}C_{lm}\biggl\{\dsize\sum_{l_1+l_2\le l,l_1\geq 1}C_{l_1 l_2}
\biggl(S^{l_1}(f_1(s,\na_x\dot\vp))\p_{z}^2
S^{l_2}\dot\vp+S^{l_1}(f_2(s,\na_x\dot\vp))
\p_{zr}^2 S^{l_2}\dot\vp\\
&\quad+S^{l_1}(f_3(s,\na_x\dot\vp))\p_{r}^2
S^{l_2}\dot\vp+S^{l_1}(\f{1}{r^2}f_4(s,\na_x\dot\vp))\p_{\th}^2
S^{l_2}\dot\vp
+S^{l_1}(\f{1}{r}f_5(s,\na_x\dot\vp))\p_{z\th}^2 S^{l_2}\dot\vp\\
&\quad +S^{l_1}(\f{1}{r}f_6(s,\na_x\dot\vp))\p_{r\th}^2
S^{l_2}\dot\vp\biggr)\biggr\}+\dsize\sum_{0\leq l\le
m}C_{lm}\biggl\{\dsize\sum_{l_1+l_2\le
l}\f{(-1)^{l_1}}{r}S^{l_2}(f_7(s,\na_x\dot\vp))\biggr\}.\endalign
$$
Let us stress that only $I_{11}$ contains derivatives of $\dot\vp$ of
order $m+2$, while $I_{12}$ and $I_{13}$ contain derivatives of
$\dot\vp$ of order at most $m+1$ which are thus lower-order terms.

>From the expressions of $f_i$, $1\leq i\leq 7$, in (3.2), the inductive
hypothesis (5.7), and Lemma 5.1, one obtains, for $m\le k_0-1$,
$$|I_{12}|+|I_{13}|\leq C(b_0,q_0)\ve\ds\sum_{0\leq l\leq m}|\na_x
S^{l}\dot\vp|^2,
$$
which implies
$$\biggl|\iiint_{D_T}(I_{12}+I_{13})\cdot\Cal M S^m\dot\vp
\,dzdrd\th\biggr|
\leq C(b_0,q_0)\ve
\iiint_{D_T}z^{\mu}\sum_{0\leq l\leq m}|\na_x S^l\dot\vp|^2
\,dzdrd\th.\tag 5.12$$

Next, we treat the troublesome term
$\ds\iiint_{D_T}I_{11}\cdot\Cal M S^m\dot\vp\,dzdrd\th$.
In view of the formulas
$\p_{y_1}^2v\p_{y_1}v=\f12\p_{y_1}((\p_{y_1}v)^2)$,
$\p_{y_1}^2v\p_{y_2}v=\p_{y_1}(\p_{y_1}v\p_{y_2}v)-\f12\p_{y_2}((\p_{y_1}v)^2)$
and $\p_{y_1y_2}^2v\p_{y_3}v$
$=\f12\{\p_{y_1}(\p_{y_2}v\p_{y_3}v)+\p_{y_2}(\p_{y_1}
v\p_{y_3}v)-\p_{y_3}(\p_{y_1}v\p_{y_2}v)\}$, one has
$$I_{11}\cdot\Cal M
S^{m}\dot\vp=\p_{z}I_{11}^1+\p_{r}I_{11}^2+\p_{\th}I_{11}^3+I_{11}^4,$$
where
$$\align
I_{11}^1&=\f{1}{2}f_1 z^{\mu}r a (\p_{z}S^m\dot\vp)^2-\f{1}{2}f_3
z^{\mu} r a (\p_r S^m\dot\vp)^2-\f{1}{2r^2}f_4 z^{\mu} r a
(\p_{\th}S^m\dot\vp)^2\\
&\qquad-\f{1}{2r}f_6 z^{\mu} ra\p_{\th}S^m\dot\vp\p_r ^m\dot\vp+f_1
z^{\mu+1}b\p_{z}S^m\dot\vp\p_{r}S^{m}\dot\vp+\f{1}{2}f_2 z^{\mu+1}
b(\p_{r}S^m\dot\vp)^2\\
&\qquad+\f{1}{2r}f_5 z^{\mu+1}
b\p_{\th}S^{m}\dot\vp\p_{r}S^{m}\dot\vp,\\
I_{11}^2&=\f{1}{2}f_2 z^{\mu} r a
(\p_{z}S^{m}\dot\vp)^2+f_3 z^{\mu} r
a\p_{r}S^{m}\dot\vp\p_{z}S^{m}\dot\vp+\f{1}{2r}f_6 z^{\mu}
r a \p_{\th}S^{m}\dot\vp\p_{z}S^{m}\dot\vp\\
&\qquad-\f{1}{2}f_1 z^{\mu+1} b(\p_{z}S^{m}\dot\vp)^2+\f{1}{2}f_3
z^{\mu+1}(\p_{r}S^{m}\dot\vp)^2-\f{1}{2r^2}f_4 z^{\mu+1}
b(\p_{\th}S^{m}\dot\vp)^2\\
&\qquad+\f{1}{2r}f_5 z^{\mu+1}
b\p_{r}S^{m}\dot\vp\p_{z}S^{m}\dot\vp,\\
I_{11}^3&=\f{1}{r}\biggl\{\f{1}{r}f_4 z^{\mu} r a
\p_{\th}S^{m}\dot\vp\p_{z}S^m\dot\vp+\f{1}{2}f_5 z^{\mu} r a
(\p_{z}S^m\dot\vp)^2+\f{1}{2}f_6 z^{\mu} r
a\p_{z}S^{m}\dot\vp\p_{r}S^{m}\dot\vp\\
&\qquad+\f{1}{r}f_4 z^{\mu+1}
b\p_{\th}S^{m}\dot\vp\p_{r}S^{m}\dot\vp+\f{1}{2}f_5 z^{\mu+1}
b\p_{z}S^{m}\dot\vp\p_{r}S^{m}\dot\vp+\f{1}{2}f_6 z^{\mu+1} b
(\p_{r}S^{m}\dot\vp)^2\biggr\},\\
I_{11}^4&=\biggl\{-\f{1}{2}\p_{z}(f_1
z^{\mu}ra)-\f{1}{2}\p_{r}(f_2 z^{\mu} r
a)-\f{1}{2}\p_{\th}(\f{1}{r}f_5 z^{\mu} r a)+\f{1}{2}\p_{z}(f_1
z^{\mu+1}b)\biggr\}(\p_{z}S^{m}\dot\vp)^2\\
&\qquad+\biggl\{-\p_{r}(f_3 z^{\mu} r a)-\f{1}{2}\p_{\th}(\f{1}{r}
z^{\mu}r a)-\p_{z}(f_1 z^{\mu+1}b)-\f{1}{2r}\p_{\th}(f_5
z^{\mu+1}b)\biggr\}\p_{z}S^{m}\dot\vp\p_{r}S^{m}\dot\vp
\endalign
$$

$$
\align
&\qquad+\biggl\{\f{1}{2}\p_{z}(f_3 z^{\mu} r a)-\f{1}{2}\p_{z}(f_2
z^{\mu+1} b)-\f{1}{2}\p_{r}(f_3 z^{\mu+1}b)-\f{1}{2r}\p_{\th}(f_6
z^{\mu+1}b)\biggr\}(\p_{r}S^{m}\dot\vp)^2\\
&\qquad+\biggl\{-\p_{\th}(\f{1}{r^2}f_4 z^{\mu}r
a)-\f{1}{2}\p_{r}(\f{1}{r}f_6 z^{\mu} r
a)+\f{1}{2}\p_{r}(\f{1}{r}f_5
z^{\mu+1}b)\biggr\}\p_{z}S^{m}\dot\vp\p_{\th}S^{m}\dot\vp\\
&\qquad+\biggl\{\f{1}{2}\p_{z}(\f{1}{r}f_6 z^{\mu}r
a)-\f{1}{r^2}\p_{\th}(f_4 z^{\mu+1}b)-\f{1}{2}\p_{z}(\f{1}{r}f_5
z^{\mu+1}b)\biggr\}\p_{r}S^{m}\dot\vp\p_{\th}S^{m}\dot\vp\\
&\qquad+\biggl\{\f{1}{2}\p_{z}(\f{1}{r^2}f_4 z^{\mu} r
a)+\p_{r}(\f{1}{r^2}f_4
z^{\mu+1}b)\biggr\}(\p_{\th}S^{m}\dot\vp)^2.
\endalign$$
Due to assumption (5.7) and in view of the properties of $f_i$, $1\leq i\leq
7$, one has
$$\ds\sum_{i=1}^{3}z^{-\mu-1}|I_{11}^i|+z^{-\mu}|I_{11}^4|\leq C(b_0,q_0)
\ve|\na_x S^{m}\dot\vp|^2,$$
which yields
$$\align
\biggl|\iiint_{D_T} & I_{11}\cdot\Cal B_0 S^{m}\dot\vp\, dzdrd\th\biggr|\\
& \leq C(b_0,q_0)\ve\biggl\{\iiint_{D_T}z^{\mu}|\na_x
S^{m}\dot\vp|^2 dzdrd\th+T^{\mu+1}\iint_{b_0 T\leq r\leq
\chi(T,\th)}|\na_x S^{m}\dot\vp|^2(T,r,\th)\, drd\th\\
& \qquad +\iint_{b_0 \leq r\leq\chi(1,\th)}|\na_x
S^{m}\dot\vp|^2(1,r,\th)\,drd\th+\iint_{\Gamma_T}z^{\mu+1}|\na_x
S^{m}\dot\vp|^2\,dS\biggr\}+\biggl|\iint_{B_T}\biggl(b_0
I_{11}^1-I_{11}^2\biggr)\,dS\biggr|.
\endalign
$$
The boundary condition $\p_{r}S^{m}\dot\vp=b_0\p_{z}S^{m}\dot\vp$ on
$B_T$, together with condition (4.6) and
$f_6|_{s=b_0}=b_0f_5|_{s=b_0}$, gives by a direct computation
$$\iint_{B_T}\biggl(b_0 I_{11}^1-I_{11}^2\biggr)dS=0,
$$
and further
$$\align
\biggl|\iiint_{D_T}I_{11} & \cdot\Cal B_0 S^{m}\dot\vp\, dzdrd\th\biggr|\\
&\leq C(b_0,q_0)\ve\biggl\{\iiint_{D_T}z^{\mu}|\na_x
S^{m}\dot\vp|^2 dzdrd\th+T^{\mu+1}\iint_{b_0 T\leq r\leq
\chi(T,\th)}|\na_x S^{m}\dot\vp|^2(T,r,\th) drd\th\\
&\qquad +\iint_{b_0 \leq r\leq\chi(1,\th)}|\na_x
S^{m}\dot\vp|^2(1,r,\th)drd\th+\iint_{\Gamma_T}z^{\mu+1}|\na_x
S^{m}\dot\vp|^2 dS\biggr\}.
\endalign$$
Combining this with (5.11), (5.12), and the initial condition (3.8)
yields
$$\align
|I_1| & \leq C(b_0,q_0)\ve\ds\sum_{0\leq l\leq
m}\biggl\{\iiint_{D_T}z^{\mu}|\na_x S^{l}\dot\vp|^2
\,dzdrd\th
 +T^{\mu+1}\iint_{b_0 T\leq r\leq
\chi(T,\th)}|\na_x S^{l}\dot\vp|^2(T,r,\th)\, drd\th \\
& \qquad +\iint_{\Gamma_T}z^{\mu+1}|\na_x S^{l}\dot\vp|^2
\,dS+\ve^2\biggr\}.\tag 5.13
\endalign$$

In particular, for $m=0$, one has
$$\align
&\biggl|\iiint_{D_T}\Cal L\dot\vp\cdot\Cal M\dot\vp\,dzdrd\th\biggr|
\leq C(b_0,q_0)\ve\biggl\{\iiint_{D_T}z^{\mu}|\na_x \dot\vp|^2
\,dzdrd\th \\
& \qquad + T^{\mu+1}\iint_{b_0 T\leq r\leq
\chi(T,\th)}|\na_x \dot\vp|^2(T,r,\th)\, drd\th
+\iint_{\Gamma_T}z^{\mu+1}|\na_x \dot\vp|^2\,dS+\ve^2\biggr\}.
\endalign
$$
Substituting this into (4.27) yields
$$ \align
&C_0T^{\mu+1}\iint_{b_0T\le
r\le\chi(T,\th)}|\na_x\dot\vp|^2(T,r,\th)\, drd\th
+C_0\iiint_{D_T}{z}^{\mu}|\na_x\dot\vp|^2\, drd\th dz\\
&\qquad+C_0\iint_{\Gamma_{T}}{z}^{\mu+1}|\na_x\dot\vp|^2\,dS \leq
C(b_0,q_0)\ve^2.\tag 5.14\endalign
$$

\vskip 0.2 true cm {\bf Step 2. Estimate of $I_2$.}
The vector fields $S_{\Gamma}$ coincide with the vector fields $S$ on
$\Gamma$ and a direct computation yields
$$\align
\biggl|[S^m_{\Gamma},\int_{0}^{1}{\hat u}_r(s_0+t\xi)dt\cdot z]\xi\biggr|
&=\biggl|\ds\sum_{m_1+m_2=m, m_1\ge
1}C_{m_1m_2}S^{m_1}\biggl(\int_{0}^{1}{\hat u}_r(s_0+t\xi)dt\cdot
z\biggr)S^{m_2}\xi\biggr|\\
&\leq C(b_0,q_0)z\biggl(M\ve\ds\sum_{1\leq l\leq
m}|S^{l}\xi|+|\xi|\biggr).
\endalign
$$
Combining this with (3.7) gives
$$\biggl|[S^m_{\Gamma},\bigl(\int_{0}^{1}{\hat
u}_r(s_0+t\xi)\,dt\bigr)z]\xi\biggr|\leq
C(b_0,q_0)z\biggl(M\ve\sum_{1\leq l\leq
m}|S^{l}\dot\vp|+|\xi|\biggr),
$$
and further
$$|I_2|\leq C(b_0,q_0)\biggl(M\ve\ds\sum_{1\leq l\leq m-1}
\iint_{\Gamma_{T}}z^{\mu+1}|\na_x S^{l}\dot\vp|^2\,ds
+\iint_{\Gamma_{T}}z^{\mu+1}|\xi|^2\,dS\biggr).\tag 5.15$$

\vskip 0.2 true cm {\bf Step 3. Estimate of $I_3$.}
One has that
$$\Cal B_0(S^{m}-S_{\Gamma}^{m})\dot\vp=[\Cal B_0, S^{m}-S_{\Gamma}^{m}]\dot\vp
+(S^{m}-S_{\Gamma}^{m})\Cal B_0\dot\vp.\tag5.16
$$
It follows from assumption (5.7) and a direct computation that
$$\biggl|[\Cal B_0, S^m-S^{m}_{\Gamma}]\dot\vp\biggr|\leq C(b_0,q_0)M\ve
\dsize\sum_{0\le l\le m-1} z^{l}|\na_x^{l+1}\dot\vp|\tag5.17
$$
and
$$|(S^{m}-S^{m}_{\Gamma})\Cal B_0\dot\vp|\leq
C(b_0,q_0)M\ve\dsize\sum_{0\leq l\leq
m}z^{l}|\na_x^{l+1}\dot\vp|.\tag5.18
$$
Combining (5.16)--(5.18) yields
$$|I_3|\leq C(b_0,q_0)M\ve\ds\sum_{0\leq l\leq
m}\iint_{\Gamma_T}z^{2l+\mu+1}|\na_x^{l+1}\dot\vp|^2\, dS.\tag
5.19$$

\vskip 0.2 true cm {\bf Step 4. Estimate of $I_4$.}
As in Step 2, one has
$$[\Cal B_0,
S^m_{\Gamma}]\dot\vp=[\Cal B_0, S^m]\dot\vp \quad \text{on $\Gamma$.}
$$
Due to
$$[\Cal B_0, S^m]\dot\vp=\dsize\sum_{0\le l\le m-1}C_{ml}S^l\Cal B_0\dot\vp$$
and
$$S^l\Cal B_0\dot\vp+\mu_2S^l\xi=S^l\kappa(\xi, \na_x\dot\vp)\quad
\text{on $\Gamma$,}
$$
one then arrives at
$$|S^l\Cal B_0\dot\vp|\le C(b_0,q_0)\biggl(|S^l\xi|+M\ve
\dsize\sum_{0\le k\le l}|S^k\xi|
+M\ve \dsize\sum_{0\le k\le l}|\na_x S^k\dot\vp|\biggr)
$$
and
$$|I_4|\le C(b_0,q_0)\biggl(\ds\sum_{0\leq l\leq m-1}\iint_{\Gamma_{T}}
z^{\mu+1}|\na_x S^{l}\dot\vp|^2 \,ds
+\iint_{\Gamma_{T}}z^{\mu+1}\xi^2 \,dS\biggr).\tag5.20$$
Substituting  (5.13)--(5.15), (5.19)--(5.20) into (5.9) and applying
Lemma 5.1 yields
$$
\align
&\iint_{b_0T\le r\le\chi(T,\th)}\dsize\sum_{0\le l\le
k_0-1}T^{2l+\mu+1} |\na_x^{l+1}\dot\vp|^2(T,r,\th)\,dS
+\iiint_{D_T}\dsize\sum_{0\le l
\le k_0-1}z^{2l+\mu}|\na_x^{l+1}\dot\vp|^2\,dzdrd\th\\
&\qquad +\iint_{\Gamma_T}\dsize\sum_{0\le l\le k_0-1}z^{2l+\mu+1}
|\na_x^{l+1}\dot\vp|^2\,dS\le
C(b_0,q_0)\biggl(\iint_{\Gamma_{T}}z^{\mu-1}|z\xi|^2
\,dS+\ve^2\biggr).\tag5.21
\endalign
$$

On the other hand, by the boundary condition (3.7), Lemma 4.2, and
(5.14), one has
$$\iint_{\Gamma_{T}}z^{\mu-1}|z\xi|^2\, dS\leq
C(b_0,q_0)\iint_{\Gamma_{T}}z^{\mu-1}\dot\vp^2 \, dS
\leq C(b_0,q_0)\ve^2.$$
Together with (5.21), this yields
$$
\align
&\iint_{b_0T\le r\le\chi(T,\th)}\dsize\sum_{0\le l\le
k_0-1}T^{2l+\mu+1} |\na_x^{l+1}\dot\vp|^2(T,r,\th)\,drd\th
+\iiint_{D_T}\dsize\sum_{0\le l
\le k_0-1}z^{2l+\mu}|\na_x^{l+1}\dot\vp|^2\,dzdrd\th\\
&\qquad +\iint_{\Gamma_T}\dsize\sum_{0\le l\le k_0-1}z^{2l+\mu+1}
|\na_x^{l+1}\dot\vp|^2\,dS\le C(b_0,q_0)\ve^2,
\endalign
$$
which shows that the conclusion of Theorem 5.3 holds. \qed

\vskip 0.2 true cm \centerline{\bf \S6. Proof of Theorem 1.1.}
\vskip 0.2 true cm

Based on the higher-order energy estimate established in Theorem 5.3,
we now prove the global existence of a conic shock wave, as asserted
in Theorem 1.1, by using a local existence result and continuous
induction. For any given $z_0>0$, the solution to (1.8) with the
initial data given on $z=z_0$ and boundary conditions (1.9)--(1.11)
exists in an interval $[z_0, z_0+\zeta]$ for some $\zeta>0$ by the
local existence result of [24] or [11, Appendix] provided that the
initial data is smooth and satisfies the compatibility
conditions. Moreover, if the perturbation of the initial data given on
$z=z_0$ is $O(\ve)$, then the lifespan of the solution is at least
$C/\ve$, with some $C>0$. Therefore, as long as one can establish that
the $L^\infty$-norms of $\dot\vp$, $\xi$, and their derivatives decay
with a rate in $z$, the solution can be continuously extended to the
whole domain.  That is, by the local existence result and the decay
properties of the solution one obtains the uniform boundedness of
$\dot\vp$, $\xi$, and their derivatives, and then one extends the
solution continuously from $z_0\le z\le z_0+\zeta$ to $z_0+\zeta\le
z\le z_0+2\zeta$, with $\zeta>0$ being independent of $z_0$. Hence,
the key to proving Theorem 1.1 is to establish the decay of the
$L^\infty$-norm of $\dot \vp$, $\xi$, and their derivatives.

It follows from Sobolev's embedding theorem (see also [11, Lemma 14])
and the assumptions of Theorem~5.4 that, for $b_0z\le r\le\chi(z,\th)$
and $1\le z\le T$, one has
$$\dsize\sum_{0\le l\le k_0-3}|z^l\na_x^{l+1}\dot\vp|^2
\le Cz^{-1}\iint_{b_0z\le r\le\chi(z,\th)} \dsize\sum_{0\le l\le
k_0-1} |z^l\na_x^{l+1}\dot\vp|^2\, drd\th.$$
On the other hand, (5.8) shows that
$$\iint_{b_0 z\le r\le\chi(z,\th)}\dsize\sum_{0\le l\le k_0-1}
|z^l\na_x^{l+1}\dot\vp|^2\,drd\th\le C(b_0,q_0)\ve^2 z^{-\mu-1}.$$
Hence, $\dsize\sum_{0\le l\le k_0-2}|z^l\na_x^{l+1}\dot\vp|^2 \le
C(b_0,q_0)\ve^2 z^{-\mu-2}$ for $b_0z\le r\le\chi(z,\th)$ and
$1\le z\le T$. For $k_0\ge 7$, one has $\dsize\sum_{l\le
[\f{k_0}{2}]+1} |z^l\na_x^{l+1}\dot\vp|\le
C(b_0,q_0)\ve{{z^{-\f{\mu}{2}-1}}}$. In addition, due to $k_0-3
\ge [\ds\f {k_0}2]+1$, Eqs.~(3.6) and (3.7) yield
$\ds\sum_{0\leq l_1+l_2\leq [\f{k_0}{2}]+1}z^{l_1}
|\p_{z}^{l_1}\p_{\th}^{l_2}\xi|\le{C(b_0,q_0)
\ve}{z^{-\f{\mu}{2}-1}}$.  If one now chooses $\mu\in (-2, -1)$, then
it follows by continuous induction that the proof of
Theorem~1.1 is complete. \qed

\vskip 0.5 true cm \centerline{\bf  Appendix} \vskip 0.3 true cm

{\bf Proof of Lemma 3.2.}  First we estimate $P_5(s)$ and $P_1'(s)$.
By (3.2), Lemmas 2.1-2.2 and Remark 2.2, then
$$\align
P_5(s)&=\f{1}{\hat u_z^2(s)-c^2(\hat \rho(s))}\biggl(-\f{\g-1}{2}s^2
\hat u_z(s)\hat u_r(s)+\f{\g+1}{2}s\hat  u_r(s) \hat u_r'(s)\\
&\qquad+s\hat u_z(s)
\hat u_z'(s)+\f{\g-1}{2}\hat u_r^2(s)-\f{1}{2}c^2(\hat \rho(s))\biggr)\\
&=\f{1+O(\c)+O(\b)}{\f{q_0^2}{(1+b_0^2)^2}-\f{\g-1}{2(1+b_0^2)}(b_0
q_0)^2}\biggl(-\f{\g-1}{2}b_0^2\cdot\f{b_0
q_0}{1+b_0^2}\cdot\f{b_0
q_0}{(1+b_0^2)^2}\\
&\qquad-\f{\g+1}{2}b_0\cdot\f{b_0
q_0}{1+b_0^2}\cdot\f{q_0}{(1+b_0^2)^2}+b_0\cdot\f{q_0}{1+b_0^2}\cdot\f{b_0
q_0}{(1+b_0^2)^2}+\f{\g-1}{2}\f{(b_0 q_0)^2}{(1+b_0^2)^2}
-\f{\g-1}{4(1+b_0^2)}(b_0 q_0)^2\biggr)\\
&=\f{-\f{\g-1}{4}b_0^2(1+b_0^2)}{1-\f{\g-1}{2}b_0^2(1+b_0^2)}
\biggl(1+O(\c)+O(\b)\biggr).
\endalign$$

Next we derive an expression for $P_1'(s)$.
Recall from (3.2) that
$$P_1(s)=\f{\hat u_{z}(s)\hat u_{r}(s)}{\hat u_z^2(s)-c^2(\hat \rho(s))}.$$
Then
$$
P_1'(s)=\f{\hat u_z'(s)\hat u_r(s)+\hat u_z(s)\hat u_r'(s)}{\hat u_z^2(s)-c^2(\hat \rho(s))}
-\f{\hat u_z(s)\hat u_r(s)(2\hat u_z(s)\hat u_z'(s)-A\g(\g-1){\hat \rho}^{\g-2}(s)
\hat \rho'(s))}{(\hat u_z^2(s)-c^2(\hat \rho(s)))^2}.
$$
It follows from Lemmas 2.1-2.2, Remark 2.2 and a direct computation that
$$\align
P_1'(s)&=\f{\f{b_0 q_0}{(1+b_0^2)^2}\f{b_0
q_0}{1+b_0^2}(1+O(\b))-\f{q_0}{(1+b_0^2)^2}
\f{q_0}{1+b_0^2}(1+O(\b))}{\f{q_0^2}{(1+b_0^2)^2}(1+O(\b))
-\f{\g-1}{2(1+b_0^2)}(b_0 q_0)^2(1+O(\c)+O(\b))}\\
& \qquad -\f{\f{q_0}{1+b_0^2}\f{b_0
q_0}{1+b_0^2}(1+O(\c)+O(\b))\biggl(2\f{q_0}{1+b_0^2}\f{b_0
q_0}{(1+b_0^2)^2}+O(\c)O(\b)O(\f{1}{b_0})\biggr)}{
\biggl(\f{q_0^2}{(1+b_0^2)^2}(1+O(\b))
-\f{\g-1}{2(1+b_0^2)}(b_0 q_0)^2(1+O(\c)+O(\b))\biggr)^2}\\
&= \f{(-1+b_0^2)\biggl(1+O(\c)+O(\b)\biggr)}{(1-\f{\g-1}{2}
b_0^2(1+b_0^2))(1+b_0^2)}
-\f{2b_0^2\biggl(1+O(\c)+O(\b)\biggr)}
{(1-\f{\g-1}{2}b_0^2(1+b_0^2))^2(1+b_0^2)}\\
&=\f{-1+\f{\g-3}{2}b_0^2-\f{\g-1}{2}b_0^6}{(1-\f{\g-1}{2}
b_0^1(1+b_0^2))^2(1+b_0^2)}\biggl(1+O(\c)+O(\b)\biggr).
\endalign$$

Other estimates in Lemma 3.2 can be carried out analogously in terms
of the expressions for $P_i$, $1\leq i\leq 5$, in (3.2) and Lemmas 2.1
and 2.2; we omit the details here. The proof of Lemma 3.2 is
complete. \qed

\smallskip

Next, we provide the proof of Lemma 3.3.

\smallskip

{\bf Proof of Lemma 3.3.} From the expressions for $B_i$, $i=1,2,3$, in
(3.5) and Lemmas 2.1 and 2.2, one has
$$\align
B_1&=-\f{(\f{\g-1}{2A\g(1+b_0^2)})^{\f{1}{\g-1}}(b_0
q_0)^{\f{2}{\g-1}}}{\f{\g-1}{2(1+b_0^2)}(b_0 q_0)^2}\f{b_0
q_0}{1+b_0^2}\biggl(\f{(b_0
q_0)^2}{(1+b_0^2)^2}-\f{q_0}{1+b_0^2}\f{b_0^2 q_0}{1+b_0^2}\biggr)
\biggl(1+O(\c)+O(\b)\biggr)\\
&\qquad+2\biggl(\f{\g-1}{2A\g(1+b_0^2)}\biggr)^{\f{1}{\g-1}}(b_0
q_0)^{\f{2}{\g-1}}\f{b_0 q_0}{1+b_0^2}(1+O(\c))\\
&=\f{2}{1+b_0^2}\biggl(\f{\g-1}{2A\g(1+b_0^2)}\biggr)^{\f{1}{\g-1}}(b_0
q_0)^{\f{\g+1}{\g-1}}\biggl(1+O(\c)+O(\b)\biggr).
\endalign$$
Similarly,
$$\align
B_2&=\biggl(\f{\g-1}{2A\g(1+b_0^2)}\biggr)^{\f{1}{\g-1}}(b_0
q_0)^{\f{2}{\g-1}}q_0\f{1-b_0^2}{1+b_0^2}\biggl(1+O(\c)+O(\b)\biggr),\\
B_3&=-\f{1}{b_0(1+b_0^2)^2}\biggl(\f{\g-1}{2A\g(1+b_0^2)}\biggr)^{\f{1}{\g-1}}
(b_0 q_0)^{\f{2\g}{\g-1}}\biggl(1+\f{1}{b_0}O(\c)+\f{1}{b_0}O(\b)\biggr).
\endalign$$
On the other hand, expressions for $\mu_i(i=1,2)$ can be
obtained from the estimates of $B_i$ which completes the proof of
Lemma 3.3. \qed

\vskip 0.3 true cm

{\bf Acknowledgments.} Yin Huicheng wishes to express his gratitude to
Professor Xin Zhouping, Chinese University of Hong Kong, and Professor
Chen Shuxing, Fudan University, Shanghai, for their constant interest
in this problem and many fruitful discussions in the past.


\bye
\enddocument